\documentclass[a4paper,12pt]{amsart}
\usepackage{amssymb,amsthm}
\usepackage[all]{xy}

\begin{document}

\newcommand{\opp}{\bowtie }
\newcommand{\po}{\text {\rm pos}}
\newcommand{\supp}{\text {\rm supp}}
\newcommand{\End}{\text {\rm End}}
\newcommand{\diag}{\text {\rm diag}}
\newcommand{\Lie}{\text {\rm Lie}}
\newcommand{\Ad}{\text {\rm Ad}}
\newcommand{\car}{\mathcal R}
\newcommand{\Tr}{\rm Tr}
\newcommand{\Spec}{\text{\rm Spec}}
\newcommand{\Cyc}{\text{\rm Cyc}}

\newtheorem{theorem}{Theorem}[section]
\newtheorem{lem}[theorem]{Lemma}
\newtheorem{cor}[theorem]{Corollary}
\newtheorem{prop}[theorem]{Proposition}
\newtheorem{thm}[theorem]{Theorem}
\newtheorem*{rmk}{Remark}
\newtheorem{eg}[theorem]{Example}

\def\ge{\geqslant}
\def\le{\leqslant}
\def\a{\alpha}
\def\b{\beta}
\def\c{\chi}
\def\g{\gamma}
\def\G{\Gamma}
\def\d{\delta}
\def\D{\Delta}
\def\L{\Lambda}
\def\e{\epsilon}
\def\et{\eta}
\def\io{\iota}
\def\o{\omega}
\def\p{\pi}
\def\ph{\phi}
\def\ps{\psi}
\def\r{\rho}
\def\s{\sigma}
\def\t{\tau}
\def\th{\theta}
\def\k{\kappa}
\def\l{\lambda}
\def\z{\zeta}
\def\v{\vartheta}
\def\x{\xi}
\def\i{^{-1}}

\def\mapright#1{\smash{\mathop{\longrightarrow}\limits^{#1}}}
\def\mapleft#1{\smash{\mathop{\longleftarrow}\limits^{#1}}}
\def\mapdown#1{\Big\downarrow\rlap{$\vcenter{\hbox{$\scriptstyle#1$}}$}}

\def\ca{\mathcal A}
\def\cb{\mathcal B}
\def\cc{\mathcal C}
\def\cd{\mathcal D}
\def\ce{\mathcal E}
\def\cf{\mathcal F}
\def\cg{\mathcal G}
\def\ch{\mathcal H}
\def\ci{\mathcal I}
\def\cj{\mathcal J}
\def\ck{\mathcal K}
\def\cl{\mathcal L}
\def\cm{\mathcal M}
\def\cn{\mathcal N}
\def\co{\mathcal O}
\def\cp{\mathcal P}
\def\cq{\mathcal Q}
\def\car{\mathcal R}
\def\cs{\mathcal S}
\def\ct{\mathcal T}
\def\cu{\mathcal U}
\def\cv{\mathcal V}
\def\cw{\mathcal W}
\def\cz{\mathcal Z}
\def\cx{\mathcal X}
\def\cy{\mathcal Y}

\def\tz{\tilde Z}
\def\tl{\tilde L}
\def\tc{\tilde C}
\def\ta{\tilde A}
\def\tx{\tilde X}

\newtheorem*{th1}{Lemma 1.3}
\newtheorem*{th2}{Lemma 1.4}
\newtheorem*{th3}{Proposition 1.7}
\newtheorem*{th4}{Corollary 1.8}
\newtheorem*{th5}{Lemma 2.2}
\newtheorem*{th6}{Proposition 2.4}
\newtheorem*{th7}{Corollary 2.5}
\newtheorem*{th8}{Corollary 2.6}
\newtheorem*{th9}{Theorem 3.2}
\newtheorem*{th10}{Proposition 3.4}
\newtheorem*{th11}{Corollary 3.5}
\newtheorem*{th12}{Lemma 3.6}
\newtheorem*{th13}{Corollary 3.7}
\newtheorem*{th14}{Corollary 3.8}
\newtheorem*{th15}{Lemma 4.1}
\newtheorem*{th16}{Proposition 4.3}
\newtheorem*{th17}{Lemma 4.4}
\newtheorem*{th18}{Corollary 4.5}
\newtheorem*{th19}{Corollary 4.6}
\newtheorem*{th20}{Proposition 5.3}
\newtheorem*{th21}{Corollary 5.4}
\newtheorem*{th22}{Corollary 5.5}
\newtheorem*{th23}{Proposition 5.6}
\newtheorem*{th24}{Lemma 5.7}
\newtheorem*{th25}{Proposition 5.8}
\newtheorem*{th26}{Corollary 5.9}
\newtheorem*{th27}{Lemma 6.3}
\newtheorem*{th28}{Corollary 6.4}
\newtheorem*{th29}{Proposition 6.5}
\newtheorem*{th30}{Proposition 6.8}
\newtheorem*{th31}{Lemma 7.2}
\newtheorem*{th32}{Lemma 7.3}
\newtheorem*{th33}{Lemma 7.4}
\newtheorem*{th34}{Theorem 7.5}
\newtheorem*{th35}{Theorem 7.6}
\newtheorem*{th36}{Corollary 7.7}
\newtheorem*{th37}{Lemma 7.9}
\newtheorem*{th38}{Lemma 7.10}
\newtheorem*{th39}{Lemma 7.11}
\newtheorem*{th40}{Lemma 7.13}
\newtheorem*{th41}{Lemma 7.15}
\newtheorem*{th42}{Lemma 7.19}
\newtheorem*{th43}{Lemma 7.21}
\newtheorem*{th44}{Lemma 7.24}
\newtheorem*{th45}{Corollary 7.25}
\newtheorem*{th46}{Theorem 7.26}

\title[Minimal length elements]
{Minimal length elements in some double cosets of Coxeter groups}
\author{Xuhua He}
\address{Department of Mathematics, Stony Brook University, Stony Brook, NY 11794, USA}%
\email{hugo@math.sunysb.edu}


\begin{abstract} We study the minimal length elements in some double cosets of Coxeter groups and use them to study Lusztig's $G$-stable pieces and the generalization of $G$-stable pieces introduced by Lu and Yakimov. We also use them to study the minimal length elements in a conjugacy class of a finite Coxeter group and prove a conjecture in \cite{GKP}.
\end{abstract}
\maketitle

\section*{Introduction}

\subsection{} Let $W$ be a Coxeter group generated by the simple reflections $s_i$ (for $i \in I$). Let $\co$ be a conjugacy class of $W$ and $\co_{\min}$ be the set of minimal length elements in $\co$. In \cite{GP1} and \cite[section 3]{GP2}, Geck and Pfeiffer obtained the following result:

If $W$ is a finite Coxeter group, then

(1) For any $w \in \co$, there exists a sequence of conjugations by $s_i$ which reduces $w$ to an element in $\co_{\min}$ and the lengths of the elements in the sequence weakly decrease;

(2) If $w, w' \in \co_{\min}$, then they are strongly conjugate in the sense of \cite[3.2.4]{GP2}.

This result was later generated by Geck, Kim and Pfeiffer to the ``twisted'' conjugacy classes of the finite Coxeter groups. See \cite{GKP}.

\subsection{} Let $J$ be a subset of $I$ and $W_J$ be the subgroup of $W$ generated by $s_j$ (for $j \in J$). The group $W_J$ acts on $W$ by conjugation. This action arises naturally in the study of Lusztig's $G$-stable pieces in \cite{L3}. A natural question is whether the above result can be generalized to the $W_J$-orbits in $W$. The answer is yes as we will see in Corollary 3.8. We will then use this result to study Lusztig's $G$-stable pieces.

\subsection{} Recently, Lu and Yakimov obtained a generalization of Lusztig's $G$-stable pieces in \cite{LY}, which is called $\car_{\cc'} \times \car_{\cc}$-stable pieces. Their motivation for studying such a generalization comes from Poisson geometry. For more details, see \cite[Introduction]{LY}.

\subsection{} In this paper, we will study these $\car_{\cc'} \times \car_{\cc}$-stable pieces in a different way. Namely, we will first study their analogy in terms of Coxeter groups. We consider the double cosets $W_{c'} \backslash (W_1 \times W_2)/W_c$, where $W_{c'}$ and $W_c$ are certain subgroups of the product $W_1 \times W_2$ of two Coxeter groups. For the minimal length elements in the double cosets, a generalization of 0.1 will be proved. Then we will use the minimal length elements to study
the $\car_{\cc'} \times \car_{\cc}$-stable pieces.

\subsection{} As an easy consequence of the results on the minimal length elements in the double cosets, we obtain some results on the minimal length elements in the (``twisted'') $W_J$-conjugacy classes on $W$, where $W_J$ is a proper parabolic subgroup of $W$. Then we will use these elements to study the (``twisted'') conjugacy classes of $W$. We will get a new proof of the results in 0.1 for finite Coxeter groups of classical type. We will also study the ``good elements'' and prove a conjecture in
\cite[5.6]{GKP}. Combining this result with the earlier results in \cite{GM} and \cite{GKP}, the existence of ``good elements'' in each twisted conjugacy class of a finite Coxeter group is established.

\subsection{} We now review the content of this paper in more detail.

In section 1, we generalize a result of B\'edard, following the approach in \cite[section 2]{L3}. In section 2, we obtain a classification of the double cosets. In section 3, we study the minimal length element in a double coset. In section 4, we introduce the notation of distinguished double cosets and distinguished elements and define a partial order on the distinguished double cosets. In section 5, we study the $\car_{\cc'} \times \car_{\cc}$-stable pieces using the distinguished elements. In section
6, we study the parabolic character sheaves and also obtain a result of the Hecke algebras. In section 7, we study the ''twisted'' conjugacy classes of finite Coxeter groups and prove the existence of the ``good elements''.

\subsection*{Acknowledgement} When I discussed with George Lusztig about the closure relations of the $G$-stable pieces, he pointed out to me the result of Geck and Pfeiffer in \cite{GP1} about the minimal length elements in conjugacy classes of finite Weyl groups. This led me to think about the possible connection between the $G$-stable piece and the minimal length length in certain subsets of the corresponding Weyl groups. It is a great pleasure to thank him here.

I would also like to thank Jiang-hua Lu for telling me her joint work \cite{LY} with Yakimov
and for helpful answers to my questions. Part of this work was done during my visit at University of Hong Kong and Hong Kong University of Science and Technology. I would like to thank them for the warm hospitality.

I would also like to thank T. A. Springer for some helpful discussions. I would also like to thank M. Geck and G. Pfeiffer for pointing out the reference \cite{F} to me.

\section{A generalization of a B\'edard's result}

In this section, we generalize a result of B\'edard \cite{Be}. We follow the approach in \cite[section 2]{L3} (and also take into account some simplification in \cite{H3}).

\subsection{} Let $I$ be a finite set and $(m_{i j})_{i, j \in I}$ be a matrix with entries in $\mathbb N \cup \{\infty\}$ such that $m_{i i}=1$ and $m_{i j}=m_{j i} \ge 2$ for all $i \neq j$. Let $W$ be a group defined by the generators $s_i$ for $i \in I$ and the relations $(s_i s_j)^{m_{i j}}=1$ for $i, j \in I$ with $m_{i j}< \infty$. We say that $(W, I)$ is a {\it Coxeter group}. Sometimes we just call $W$ itself a Coxeter group.

We denote by $l$ the length function and $\le$ the Bruhat order. For $w \in W$, we denote by $\supp(w)$ the set of simple reflections that appear in some (or equivalently, any) reduced expression of $w$. For $J \subset I$, we denote by $W_J$ the standard parabolic subgroup of $W$ generated by $J$ and by $W^J$ (resp. ${}^J W$) the set of minimal coset representatives in $W/W_J$ (resp. $W_J \backslash W$). For $J, K \subset I$, we simply write $W^J \cap {}^K W$ as ${}^K W^J$.

For $J \subset I$ with $W_J$ finite, we denote by $w_J$ the maximal element in $W_J$.

For an automorphism $\s$ of $W$, we define the $\s$-twisted conjugation action of $W$ on itself by $w \cdot w'=w w' \s(w) \i$. The orbits are called the $\s$-twisted conjugacy classes of $W$.

For a finite set $X$, we denote by $\sharp X$ its cardinality.

\subsection{} We recall some known results about
$W^J$.

(1) If $w \in W^J$ and $i \in I$, then there are three
possibilities.

$\quad$ (a) $s_i w>w$ and $s_i w \in W^J$;

$\quad$ (b) $s_i w>w$ and $s_i w=w s_j$ for some $j \in J$;

$\quad$ (c) $s_i w<w$ in which case $s_i w \in W^J$.

(2) If $w \in W^J$, $v \in W_J$ and $K \subset J$, then $v \in
W^K$ if and only if $w v \in W^K$.

(3) If $w \in {}^{J'} W^J$ and $u \in W_{J'}$, then $u w \in W^J$
if and only if $u \in W^K$, where $K=J' \cap \Ad(w) J$.

\begin{th1} (1) Let $J, K \subset I$ and $w \in {}^K W$ with $w \i(K) \subset J$. Assume that $w=x y$ for $x \in {}^K W^J$ and $y \in W_J$. Then $x \i(K) \subset J$.

(2) Let $J, K \subset I$ and $w \in W^K$ with $w(K) \subset J$. Assume that $w=x y$ for $x \in W_J$ and $y \in {}^J W^K$. Then $y(K) \subset J$.
\end{th1}

We only prove part (1). Part (2) can be proved in the same way.

By assumption, for $k \in K$, there exists $j \in J$, such that $s_k w=w s_j=x y s_j$. It is easy to see that $s_k x>x$. If $s_k x \in W^J$, then $$s_k w \in (s_k x) W_J, x y s_j \in x W_J \text{ and } s_k x, x \in W^J.$$ Thus $s_k x=x$, which is a contradiction. Hence by 1.2 (1), $s_k x=x s_{j'}$ for some $j' \in J$. The lemma is proved.

\begin{th2} Let $u, w \in W$. Then

(1) The subset $\{v w; v \le u\}$ of $W$ contains a unique
minimal element $y$. Moreover, $l(y)=l(w)-l(y w \i)$.

(2) The subset $\{v w; v \le u\}$ of $W$ contains a unique
maximal element $y'$. Moreover, $l(y')=l(w)+l(y' w \i)$.
\end{th2}

\begin{rmk} This is a generalization of the result \cite[Lemma 3.3]{H1}. In {\it loc. cit.} the Coxeter group $W$ is a finite Weyl group. But this assumption is not needed here.
\end{rmk}

We will only prove part (1). Part (2) can be proved in the same way.

We argue by induction on $l(u)$. For $l(u)=0$, part (1) is clear. Assume now that $l(u)>0$ and that the statement holds for all $u' \in W$ with $l(u')<l(u)$. Then there exists $i \in I$ such that $s_i u<u$. We denote $s_i u$ by $u'$. Then by induction hypothesis, the subset $\{v' w; v' \le u'\}$ contains a unique element $y_1=v_1 w$ and $l(y_1)=l(w)-l(v_1)$. Set $y=\min\{y_1, s_i y_1\}$. Then we have that $y<s_i y$ and $y \le y_1$. Now assume that $z$ is an element in $\{v w; v \le u\}$. Then it is easy
to see that either $z$ or $s_i z$ is contained in $\{v' w; v \le u'\}$. Therefore we have that either $y \le y_1 \le z$ or $y \le y_1 \le s_i z$. In the second case, by \cite[Corollary 2.5]{L1}, we still have that $y \le z$.

So $y$ is the minimal element in $\{v w; v \le u\}$.

If $y=y_1$, then $l(y)=l(w)-l(v_1)$. If $y=s_i y_1=s_i v_1 w$, then $l(y)=l(y_1)-1=l(w)-l(v_1)-1$. Since $l(y) \ge l(w)-l(s_i v_1)$, we have that $l(s_i v_1)=l(v_1)+1$ and $l(y)=l(w)-l(s_i v_1)$. The lemma is proved.

\subsection*{1.5} Let $(W_1, I_1)$ and $(W_2, I_2)$ be two Coxeter groups.
A triple $c=(J_1, J_2, \d)$ consisting of $J_1 \subset
I_1$, $J_2 \subset I_2$ and an isomorphism $\d: W_{J_1} \rightarrow
W_{J_2}$ which sends $J_1$ to $J_2$ will be called {\it an
admissible triple for $W_1 \times W_2$}. To each admissible triple
$c=(J_1, J_2, \d)$, set $$W_c=\{(w, \d(w)); w \in W_{J_1}\} \subset
W_1 \times W_2.$$

Let $c=(J_1, J_2, \d)$ be an admissible triple for $W_1 \times W_2$,
then $c \i=(J_2, J_1, \d \i)$ is an admissible triple for $W_2
\times W_1$.

For admissible triples $c=(J_1, J_2, \d)$ and $c'=(J'_1, J'_2, \d')$
for $W_1 \times W_2$, we say that $c' \le c$ if $J'_1 \subset J_1$,
$J'_2 \subset J_2$ and $\d'=\d \mid_{W_{J'_1}}$.

\subsection*{1.6} Let $c=(J_1, J_2, \d)$ and $c'=(J'_1, J'_2, \d')$ be two admissible triples.
Let $\ct(c, c')$ be the set of all sequences $(J_1^{(n)}, J_2^{'
(n)}, w_1^{(n)}, w_2^{(n)})_{n \ge 0}$ where $J_1^{(n)} \subset
J_1$, $J_2^{' (n)} \subset J'_2$, $w_1^{(n)} \in W_1$ and $w_2^{(n)}
\in W_2$ are such that

\begin{gather*} \tag{a} J_1^{(0)}=J_1, J_2^{' (0)}=\d' w_1^{(0)} J_1 \cap J'_2; \\ \tag{b} J_1^{(n)}=\d \i (w_2^{(n-1)}) \i J_2^{' (n-1)} \cap J_1 \text{ for } n \ge 1; \\ \tag{c} J_2^{' (n)}=\d' w_1^{(n)} J_1^{(n)} \cap J'_2 \text{ for } n \ge 1; \\ \tag{d} w_1^{(n)} \in {}^{J'_1} W_1^{J_1^{(n)}}, w_2^{(n)} \in {}^{J_2^{' (n)}} W_2^{J_2} \text{ for } n \ge 0; \\ \tag{e} w_1^{(n)} \in w_1^{(n-1)} W_{J_1^{(n-1)}}, w_2^{(n)} \in W_{J_2^{' (n-1)}} w_2^{(n-1)} \text{ for } n \ge 1.
\end{gather*}

\begin{th3} $(J_1^{(n)}, J_2^{' (n)}, w_1^{(n)}, w_2^{(n)})_{n \ge 0} \mapsto (w_1^{(m)}, w_2^{(m)})$ for $m \gg 0$ is a well-defined bijection $\phi: \ct(c, c') \rightarrow {}^{J'_1} W_1 \times W_2^{J_2}$.
\end{th3}

Let $(J_1^{(n)}, J_2^{' (n)}, w_1^{(n)}, w_2^{(n)})_{n \ge 0} \in \ct(c, c')$. We prove by induction on $n \ge 0$ that
\begin{equation*} \tag{a} J_1^{(n+1)} \subset J_1^{(n)}, J_2^{' (n+1)} \subset J_2^{' (n)}.
\end{equation*}

For $n=0$, $J_1^{(1)} \subset J_1^{(0)}=J_1$. Now $(w_1^{(1)}) \i (\d') \i J_2^{' (1)} \subset J_1^{(1)} \subset J_1$ and $w_1^{(0)}=\min(w_1^{(1)} W_{J_1^{(0)}})$. By Lemma 1.3, $(w_1^{(0)}) \i (\d') \i J_2^{' (1)} \subset J_1$. Hence $J_2^{' (1)} \subset \d' w_1^{(0)} J_1 \cap J'_2=J_2^{' (0)}$.

Assume now that $n>0$ and that (a) holds when $n$ is replaced by $n-1$. Then $w_2^{(n)} \d J_1^{(n+1)} \subset J_2^{' (n)} \subset J_2^{' (n-1)}$ and $w_2^{(n-1)}=\min(W_{J_2^{' (n-1)}} w_2^{(n)})$. By Lemma 1.3, $w_2^{(n-1)} \d J_1^{(n+1)} \subset J_2^{' (n-1)}$. Hence $$J_1^{(n+1)} \subset \d \i (w_2^{(n-1)}) \i J_2^{' (n-1)} \cap J_1=J_1^{(n)}.$$

Similarly, $J_2^{' (n+1)} \subset J_2^{' (n)}$.

(a) is proved.

Now since $I_1, I_2$ are finite sets, there exists $n_0 \ge 1$ such that $J_1^{(n)}=J_1^{(n-1)}$ and $J_2^{' (n)}=J_2^{' (n-1)}$ for $n \ge n_0$. For such $n$ we have $$w_1^{(n)} \in {}^{J'_1} W_1^{J_1^{(n)}}, w_1^{(n-1)} \in {}^{J'_1} W_1^{J_1^{(n)}}, w_1^{(n)} \in w_1^{(n-1)} W_{J_1^{(n)}}.$$ Thus $w_1^{(n)}=w_1^{(n-1)}$. Similarly $w_2^{(n)}=w_2^{(n-1)}$. Thus $\phi$ is well-defined. We set $w_1=w_1^{(m)}$ and $w_2=w_2^{(m)}$ for $m \gg 0$.

By 1.6 (a) and (d), $w_1 \in w_1^{(n)} W_{J_1^{(n)}}$. Since $w_1^{(n)} \in W^{J_n^{(n)}}$, we have that \begin{equation*} \tag{b} w_1^{(n)}=\min(w_1 W_{J_1^{(n)}}).\end{equation*} Similarly, \begin{equation*}\tag{c} w_2^{(n)}=\min(W_{J_2^{' (n)}} w_2).\end{equation*}

Now assume that $\phi \bigl( (\tilde J_1^{(n)}, \tilde J_2^{' (n)}, \tilde w_1^{(n)}, \tilde w_2^{(n)})_{n \ge 0} \bigr)=(w_1, w_2)$. We show by induction on $n \ge 0$ that \begin{equation*} \tag{d} J_1^{(n)}=\tilde J_1^{(n)}, J_2^{' (n)}=\tilde J_2^{' (n)}, w_1^{(n)}=\tilde w_1^{(n)}, w_2^{(n)}=\tilde w_2^{(n)}.
\end{equation*}

For $n=0$ this holds since \begin{gather*} J_1^{(0)}=\tilde J_1^{(0)}=J_1, w_1^{(0)}=\tilde w_1^{(0)}=\min(w_1 W_{J_1}), \\ J_2^{' (0)}=\tilde J_2^{' (0)}=\d' w_1^{(0)} J_1 \cap J'_2, w_2^{(0)}=\tilde w_2^{(0)}=\min(W_{J_2^{' (0)}} w_2).\end{gather*}

Assume now that $n>0$ and that (d) holds when $n$ is replaced by $n-1$. From 1.6 (b), we deduce that $J_1^{(n)}=\tilde J_1^{(n)}=\d \i (w_2^{(n-1)}) \i J_2^{' (n-1)} \cap J_1$. From (b), we deduce that $w_1^{(n)}=\tilde w_1^{(n)}=\min(w_1 W_{J_1^{(n)}})$. By 1.6 (c), we deduce that $J_2^{' (n)}=\tilde J_2^{' (n)}$. From (c), we deduce that $w_2^{(n)}=\tilde w_2^{(n)}$.

Thus (d) holds and $\phi$ is injective.

We define an inverse to $\phi$. Let $(w_1, w_2) \in {}^{J'_1} W_1 \times W_2^{J_2}$, we define by induction on $n \ge 0$ a sequence $(J_1^{(n)}, J_2^{' (n)}, w_1^{(n)}, w_2^{(n)})_{n \ge 0}$ as follows.

We set $J_1^{(0)}=J_1$, $w_1^{(0)}=\min(w_1 W_{J_1})$, $J_2^{' (0)}=\d' w_1^{(0)} J_1 \cap J'_2$ and $w_2^{(0)}=\min(W_{J_2^{' (0)}} w_2)$.

Assume now that $n>0$ and $J_1^{(n-1)}, J_2^{' (n-1)}, w_1^{(n-1)}, w_2^{(n-1)}$ are defined. We define $J_1^{(n)}=\d \i (w_2^{(n-1)}) \i J_2^{' (n-1)} \cap J_1$, $w_1^{(n)}=\min(w_1 W_{J_1^{(n)}})$, $J_2^{' (n)}=\d' w_1^{(n)} J_1^{(n)} \cap J'_2$ and $w_2^{(n)}=\min(W_{J_2^{' (n)}} w_2)$. This completes the inductive definition.

Now for $n \ge 1$, $w_1^{(n)} \in w_1 W_{J_1^{(n)}}$ and $w_1 \in w_1^{(n-1)} W_{J_1^{(n-1)}}$. Hence $w_1^{(n)} \in \bigl(w_1^{(n-1)} W_{J_1^{(n-1}}) \bigr) W_{J_1^{(n)}}=w_1^{(n-1)} W_{J_1^{(n-1)}}$.

Similarly, $w_2^{(n)} \in W_{J_2^{' (n-1)}} w_2^{(n-1)}$.

For $n \ge 0$, $w_1=w_1^{(n)} x$ for some $x \in W_{J_1^{(n)}}$ and $l(w_1)=l(w_1^{(n)})+l(x)$. Now for $v \in W_{J'_1}$, $l(v w_1)=l(v)+l(w_1)$ since $w_1 \in {}^{J'_1} W_1$. On the other hand, $l(v w_1^{(n)} x) \le l(v w_1^{(n)})+l(x)$. Then $l(v)+l(w_1^{(n)})=l(v w_1^{(n)})$ and $w_1^{(n)} \in {}^{J'_1} W_1$.

Similarly, $w_2^{(n)} \in W_2^{J_2}$.

Thus $(J_1^{(n)}, J_2^{' (n)}, w_1^{(n)}, w_2^{(n)})_{n \ge 0} \in \ct(c, c')$.

We show that $w_1^{(m)}=w_1$ and $w_2^{(m)}=w_2$ for $m \gg 0$.

For any $n \ge 0$, we have $w_1=w_1^{(n)} u$ and $w_2=v w_2^{(n)}$ for $u \in W_{J_1^{(n)}}$ and $v \in W_{J_2^{' (n)}}$. Since $w_1 \in {}^{J'_1} W_1$ and $w_1^{(n)} \in {}^{J'_1} W_1^{J_1^{(n)}}$, by \cite[2.1(b)]{L3}, $u \in {}^{J_1^{(n)} \cap (w_1^{(n)}) \i J'_1} W_1$.

Assume that $n \gg 0$, we have $J_1^{(n)}=J_1^{(n-1)}$ and $J_2^{' (n)}=J_2^{' (n-1)}$. By 1.6 (b) and (c), \begin{gather*} \sharp J_1^{(n)} \le \sharp \bigl( w_2 \i J_2^{' (n)} \cap J_2 \bigr) \le \sharp J_2^{' (n)}, \\ \sharp J_2^{' (n)} \le \sharp \bigl( w_1 J_1^{(n)} \cap J'_1 \bigr) \le \sharp J_1^{(n)}. \end{gather*}

Hence $\sharp J_1^{(n)}=\sharp J_2^{' (n)}$ and \begin{gather*} \tag{e} w_1 J_1^{(n)} \subset J'_1, w_2 \i J_2^{' (n)} \subset J_2, \\ \tag{f} J_1^{(n)}=\d \i w_2 \i J_2^{' (n)}, J_2^{' (n)}=\d' w_1 J_1^{(n)}.\end{gather*}
So $u \in {}^{J_1^{(n)} \cap w_1 \i J'_1} W_1={}^{J_1^{(n)}} W_1$. Notice that $u \in W_{J_1^{(n)}}$. Thus $u=1$ and $w_1=w_1^{(n)}$. Similarly, $w_2=w_2^{(n)}$.

Thus we have defined a map $\psi: {}^{J'_1} W_1 \times W_2^{J_2} \rightarrow \ct(c, c')$ such that $\phi \circ \psi=id$. Hence $\phi$ is bijective. The proposition is proved.

\begin{th4} For $w_1 \in {}^{J'_1} W_1$ and $w_2 \in W_2^{J_2}$, define $$I(w_1, w_2, c, c')=\max\{K \subset J_1; w_1(K) \subset J'_1 \text{ and } \d' w_1 K=w_2 \d K\}.$$ Let
$(J_1^{(n)}, J_2^{' (n)}, w_1^{(n)}, w_2^{(n)})_{n \ge 0}$ be the element in $\ct(c, c')$ whose image under $\phi$ is $(w_1, w_2)$. Then $I(w_1, w_2, c, c')=J_1^{(n)}$ for $n \gg 0$.
\end{th4}

By (e) and (f) in the proof of the proposition 1.7, we have that $J_n^{(n)} \subset J_1, w_1(J_1^{(n)}) \subset J'_1$ and $\d' w_1 J_1^{(n)}=w_2 \d J_1^{(n)}$ for $n \gg 0$. Thus $J_1^{(n)} \subset I(w_1, w_2, c, c')$ for $n \gg 0$.

Now set $I'(w_1, w_2, c, c')=\d' w_1 I(w_1, w_2, c, c')$. It suffices to prove that for  $n \ge 0$
\begin{equation*}
I(w_1, w_2, c, c') \subset J_1^{(n)} \text{ and } I'(w_1, w_2, c, c') \subset J_2^{' (n)}.
\end{equation*}

We argue by induction on $n$. For $n=0$, $I(w_1, w_2, c, c') \subset J_1^{(0)}=J_1$. Now $$w_1 I(w_1, w_2, c, c') \subset J'_1, w_1 \i w_1 I(w_1, w_2, c, c')=I(w_1, w_2, c, c') \subset J_1.$$ Notice that $w_1^{(0)}=\min(w_1 W_{J_1})$. By Lemma 1.3, $$(w_1^{(0)}) \i w_1 I(w_1, w_2, c, c') \subset J_1.$$ In other words, $w_1 I(w_1, w_2, c, c') \subset w_1^{(0)} J_1 \cap J'_1$ and $I'(w_1, w_2, c, c') \subset \d' w_1^{(0)} J_1 \cap J'_2$. So (a) holds for $n=0$.

Assume now that $n>0$ and that (a) holds when $n$ is replaced by $n-1$. Then $w_2 \d I(w_1, w_2, c, c')=I'(w_1, w_2, c, c') \subset J_2^{' (n-1)}$ and $w_2^{(n-1)}=\min(W_{J_2^{'(n-1)}} w_2)$. By Lemma 1.3, $w_2^{(n-1)} \d I(w_1, w_2, c, c') \subset J_2^{' (n-1)}$. Hence $I(w_1, w_2, c, c') \subset \d \i (w_2^{(n-1)}) \i J_2^{' (n-1)} \cap J_1=J_1^{(n)}$.

Notice that $w_1 \i w_1 I(w_1, w_2, c, c')=I(w_1, w_2, c, c') \subset J_1^{(n)}$ and $w_1^{(n)} \in \min(w_1 W_{J_1^{(n)}})$. By Lemma 1.3, $(w_1^{(n)}) \i w_1 I(w_1, w_2, c, c') \subset J_1^{(n)}$. Hence $I'(w_1, w_2, c, c')=\d' w_1 I(w_1, w_2, c, c') \subset \d' w_1^{(n)} J_1^{(n)} \cap J'_2=J_2^{' (n)}$.

The Corollary is proved.

\subsection*{1.9} Below is a variant of the above results.

Let $\ct'(c, c')$ be the set of sequences $(J_1^{(n)}, J_2^{' (n)}, w_1^{(n)}, w_2^{(n)})_{n \ge 0}$ where $J_1^{(n)} \subset J_1$, $J_2^{' (n)} \subset J'_2$, $w_1^{(n)} \in W_1$ and $w_2^{(n)} \in W_2$ are such that

\begin{gather*} \tag{a} J_1^{(0)}=\d \i (w_1^{(0)}) \i J'_2 \cap J_1, J_2^{' (0)}=J'_2; \\ \tag{b} J_1^{(n)}=\d \i (w_2^{(n)}) \i J_2^{' (n)} \cap J_1 \text{ for } n \ge 1, \\ \tag{c} J_2^{' (n)}=\d' w_1^{(n-1)} J_1^{(n-1)} \cap J'_2 \text{ for } n \ge 1; \\ \tag{d} w_1^{(n)} \in {}^{J'_1} W_1^{J_1^{(n)}}, w_2^{(n)} \in {}^{J_2^{' (n)}} W_2^{J_2} \text{ for } n \ge 0; \\ \tag{e} w_1^{(n)} \in w_1^{(n-1)} W_{J_1^{(n-1)}}, w_2^{(n)} \in W_{J_2^{' (n-1)}} w_2^{(n-1)} \text{ for } n \ge 1.
\end{gather*}

Then $(J_1^{(n)}, J_2^{' (n)}, w_1^{(n)}, w_2^{(n)})_{n \ge 0} \mapsto (w_1^{(m)}, w_2^{(m)})$ for $m \gg 0$ is also a well-defined bijection $\ct(c, c') \rightarrow {}^{J'_1} W_1 \times W_2^{J_2}$ and $I(w_1, w_2, c, c')=J_1^{(n)}$ for $n \gg 0$.

\section{The $W_{c'} \times W_c$-stable pieces}

\subsection{} To each element $(w_1, w_2) \in W_1 \times W_2$ we associate a sequence $(J_1^{(n)}, J_2^{' (n)}, w_1^{(n)}, w_2^{(n)}, u_1^{(n)}, u_2^{(n)}, v_1^{(n)}, v_2^{(n)})_{n \ge 0}$ with $J_1^{(n)} \subset J_1$, $J_2^{' (n)} \subset J'_2$, $w_1^{(n)} \in {}^{J'_1} W_1^{J_1^{(n)}}$, $w_2^{(n)} \in {}^{J_2^{' (n)}} W_2^{J_2}$, $u_1^{(n)} \in {}^{J'_1} W_1$, $u_2^{(n)} \in W_2^{J_2}$, $v_1^{(n)} \in W_1$ and $v_2^{(n)} \in W_2$. We set
\begin{gather*}
J_1^{(0)}=J_1, u_1^{(0)}=\min(W_{J'_1} w_1), w_1^{(0)}=\min(u_1^{(0)} W_{J_1}), v_1^{(0)}=w_1 (u_1^{(0)}) \i, \\ J_2^{' (0)}= \d' w_1^{(0)} J_1 \cap J'_2, u_2^{(0)}=\min(\d'(v_1^{(0)}) \i w_2 W_{J_2}), \\ w_2^{(0)}=\min(W_{J_2^{' (0)}} u_2^{(0)}), v_2^{(0)}=(u_2^{(0)}) \i \d'(v_1^{(0)}) \i w_2.
\end{gather*}

Assume that $n>1$ and that $J_1^{(n-1)}, J_2^{' (n-1)}, w_1^{(n-1)}, w_2^{(n-1)}, u_1^{(n-1)}$, $u_2^{(n-1)}, v_1^{(n-1)}, v_2^{(n-1)}$ are already defined. Let
\begin{gather*}
J_1^{(n)}=\d \i (w_2^{(n-1)}) \i J_2^{' (n-1)} \cap J_1, u_1^{(n)}=\min\bigl(W_{J'_1} u_1^{(n-1)} \d \i (v_2^{(n-1)}) \i \bigr), \\ w_1^{(n)}=\min(u_1^{(n)} W_{J_1^{(n)}}), v_1^{(n)}=u_1^{(n-1)} \d \i(v_2^{(n-1)}) \i (u_1^{(n)}) \i, \\ J_2^{' (n)}=\d' w_1^{(n)} J_1^{(n)} \cap J'_2, u_2^{(n)}=\min \bigl( \d'(v_1^{(n)}) \i u_2^{(n-1)} W_{J_2} \bigr), \\ w_2^{(n)}=\min(W_{J_2^{' (n)}} u_2^{(n)}), v_2^{(n)}=(u_2^{(n)}) \i \d'(v_1^{(n)}) \i u_2^{(n-1)}.
\end{gather*}

This completes the inductive definition.

\begin{th5} We keep the notation in 2.1. Then $$(J_1^{(n)}, J_2^{' (n)}, w_1^{(n)}, w_2^{(n)})_{n \ge 0} \in \ct(c, c').$$
\end{th5}

1.6 (a), (b), (c) and (d) are automatically satisfied. By definition, $v_1^{(0)} \in W_{J'_1}$ and $v_2^{(0)} \in W_{J_2}$. Now we prove by induction on $n>1$ that \begin{gather*}
\tag{a} v_1^{(n)} \in W_{(\d') \i J_2^{' (n-1)}}, v_2^{(n)} \in W_{\d J_1^{(n)}}, \\ \tag{b} w_1^{(n)} \in w_1^{(n-1)} W_{J_1^{(n-1)}}, w_2^{(n)} \in W_{J_2^{' (n-1)}} w_2^{(n-1)} \end{gather*}

For $n=1$, we have \begin{gather*} u_1^{(1)}=\min\bigl(W_{J'_1} u_1^{(0)} \d \i(v_2^{(0)}) \i \bigr), \\ u_1^{(0)} \d \i(v_2^{(0)}) \i \in u_1^{(0)} W_{J_1^{(0)}}=w_1^{(0)} W_{J_1^{(0)}}.\end{gather*} Notice that $w_1^{(0)} \in {}^{J'_1} W_1^{J_1^{(0)}}$.
By \cite[Lemma 3.6]{H1}, $$u_1^{(0)} \d \i(v_2^{(0)}) \i \in W_{J'_1 \cap w_1^{(0)} J_1^{(0)}} u_1^{(1)}=W_{(\d') \i J_2^{' (0)}} u_1^{(1)}, u_1^{(1)} \in w_1^{(0)} W_{J_1^{(0)}}.$$ Hence $v_1^{(1)} \in W_{(\d') \i J_2^{' (0)}}$ and $w_1^{(1)} \in w_1^{(0)} W_{J_1^{(0)}}$.

We also have that \begin{gather*} u_2^{(1)}=\min \bigl( \d'(v_1^{(1)}) \i u_2^{(0)} W_{J_2} \bigr), \\ \d'(v_1^{(1)}) \i u_2^{(0)} \in W_{J_2^{' (0)}} u_2^{(0)}=W_{J_2^{' (0)}} w_2^{(0)}.\end{gather*} Notice that $w_2^{(0)} \in {}^{J_2^{' (0)}} W_2^{J_2}$. By \cite[Lemma 3.6]{H1}, $$\d'(v_1^{(1)}) \i u_2^{(0)} \in u_2^{(1)} W_{J_2 \cap (w_2^{(0)}) \i J_2^{' (0)}}=u_2^{(1)} W_{\d J_1^{(1)}}, u_2^{(1)} \in W_{J_2^{' (0)}} w_2^{(0)}.$$ Hence $v_2^{(1)} \in W_{\d J_1^{(1)} }$ and $w_2^{(1)} \in W_{J_2^{'
(0)}} w_2^{(0)}$.

Assume now that $n>2$ and that (a) and (b) hold when $n$ is replaced by $n-1$. Then we can show in the same way that $v_1^{(n)} \in W_{(\d') \i J_2^{' (n-1)}}$, $v_2^{(n)} \in W_{\d J_2^{(n)}}$, $w_1^{(n)} \in w_1^{(n-1)} W_{J_1^{(n-1)}}$ and $w_2^{(n)} \in W_{J_2^{' (n-1)}} w_2^{(n-1)}$. The lemma is proved.

\subsection*{2.3} We define a map $\pi: W_1 \times W_2 \rightarrow {}^{J'_1} W_1 \times W_2^{J_2}$ as follows.

Let $(J_1^{(n)}, J_2^{' (n)}, w_1^{(n)}, w_2^{(n)}, u_1^{(n)}, u_2^{(n)}, v_1^{(n)}, v_2^{(n)})_{n \ge 0}$ be the sequence associated to $(w_1, w_2) \in W_1 \times W_2$. By the previous lemma, $$(J_1^{(n)}, J_2^{' (n)}, w_1^{(n)}, w_2^{(n)})_{n \ge 0} \in \ct(c, c').$$  Now set $$\pi(w_1, w_2)=\phi \bigl((J_1^{(n)}, J_2^{' (n)}, w_1^{(n)}, w_2^{(n)})_{n \ge 0}  \bigr).$$ This completes the definition.

For $(w_1, w_2) \in {}^{J'_1} W_1 \times W_2^{J_2}$, set $[w_1, w_2, c, c']=\pi \i(w_1, w_2)$. Then $$W_1 \times W_2=\sqcup_{(w_1, w_2) \in {}^{J'_1} W_1 \times W_2^{J_2}} [w_1, w_2, c, c'].$$

\begin{th6} Let $(w_1, w_2) \in {}^{J'_1} W_1 \times W_2^{J_2}$. Then

(1) $[w_1, w_2, c, c']=W_{c'} (w_1 W_{I(w_1, w_2, c, c')}, w_2) W_c$.

(2) Define an automorphism $\s: W_{I(w_1, w_2, c, c')} \rightarrow
W_{I(w_1, w_2, c, c')}$ by $$\s(w)=\d \i \bigl(w_2 \i \d'(w_1 w w_1 \i) w_2 \bigr).$$ Then map $W_{I(w_1, w_2, c, c')} \rightarrow W_1
\times W_2$ defined by $w \rightarrow (w_1 w, w_2)$ induces a
bijection between the $\s$-twisted conjugacy classes on $W_{I(w_1,
w_2, c, c')}$ and the double cosets $W_{c'} \backslash [w_1, w_2, c,
c'] /W_c$.
\end{th6}

\begin{rmk} By part (1), for each $(w_1, w_2) \in {}^{J'_1} W_1 \times W_2^{J_2}$, the subset $[w_1, w_2, c, c']$ of $W_1 \times W_2$ is stable under the action of $W_{c'} \times W_c$. We call $[w_1, w_2, c, c']$ a $W_{c'} \times W_c$-stable piece of $W_1 \times W_2$.
\end{rmk}

Let $(J_1^{(n)}, J_2^{' (n)}, w_1^{(n)}, w_2^{(n)}, u_1^{(n)}, u_2^{(n)}, v_1^{(n)}, v_2^{(n)})_{n \ge 0}$ be the sequence associated to $(w'_1, w'_2) \in W_1 \times W_2$.

By definition, $u_1^{(n)} \in w_1^{(n)} W_{J_1^{(n)}}$ for $n \ge 0$. By (e) in the proof of Proposition 1.7, $w_1^{(n)} W_{J_1^{(n)}} \subset W_{J'_1} w_1^{(n)}$ for $n \gg 0$. Since $u_1^{(n)}, w_1^{(n)} \in {}^{J'_1} W_1$, we have that \begin{equation*}
\tag{a} u_1^{(n)}=w_1^{(n)} \text{ for } n \gg 0. \end{equation*}

Similarly, \begin{equation*}
\tag{b} u_2^{(n)}=w_2^{(n)} \text{ for } n \gg 0. \end{equation*}

By definition, $(w'_1, w'_2)$ and $(u_1^{(0)}, \d'(v_1^{(0)}) \i w'_2)$ are in the same $W_{c'} \times W_c$-coset. We can show by induction on $n \ge 0$ that \begin{align*}
\tag{c} & (w'_1, w'_2), (u_1^{(n)} \d \i(v_2^{(n)}) \i, u_2^{(n)}) \text{ and } (u_1^{(n+1)}, \d'(v_1^{(n+1)}) \i u_2^{(n)}) \\ & \text{ are in the same } W_{c'} \times W_c \text{-coset}.
\end{align*}

Now suppose that $\pi(w'_1, w'_2)=(w_1, w_2)$. For $n \gg 0$, we have $u_1^{(n)}=w_1$ and $u_2^{(n)}=w_2$. Moreover, $\d \i (v_2^{(n)}) \in W_{J_1^{(n)}}=W_{I(w_1, w_2, c, c')}$. Thus by (c), $(w'_1, w'_2) \in [w_1, w_2, c, c']$. On the other hand, it is easy to see that for $v \in W_{I(w_1, w_2, c, c')}$, $\pi(w_1 v, w_2)=(w_1, w_2)$. Therefore $[w_1, w_2, c, c']=W_{c'} (w_1 W_{I(w_1, w_2, c, c')}, w_2) W_c$.

Part (1) is proved.

For all $x \in W_{w_1 I(w_1, w_2, c, c')}$ and $y \in W_{I(w_1, w_2, c, c')}$, $$(x, \d'(x)) (w_1 W_{I(w_1, w_2, c, c')}, w_2) (y, \d(y))=(w_1 W_{I(w_1, w_2, c, c')}, w_2 W_{\d I(w_1, w_2, c, c')}).$$ On the other hand, assume that $(x, \d'(x)) (w_1 v, w_2) (y, \d(y))=(w_1 v', w_2)$ for $x \in W_{J'_1}$, $y \in W_{J_1}$ and $v, v' \in W_{I(w_1, w_2, c, c')}$, then $x w_1 v y=w_1 v'$ and $\d'(x) w_2 \d(y)=w_2$. By \cite[Lemma 3.6]{H1}, $\d' \supp(x)=w_2 \d \supp(y)$ and $$w_1 \i \bigl(\supp(x) \cup w_1
I(w_1, w_2, c, c') \bigr)=\supp(y) \cup I(w_1, w_2, c, c').$$ Therefore, \begin{gather*} w_1 \bigl(\supp(y) \cup I(w_1, w_2, c, c') \bigr) \in J'_1, \\ \d' w_1  \bigl(\supp(y) \cup I(w_1, w_2, c, c') \bigr)=w_2 \d  \bigl(\supp(y) \cup I(w_1, w_2, c, c') \bigr).\end{gather*} Hence $\supp(y) \subset I(w_1, w_2, c, c')$ and $\d'(x) w_2 \d(y)=w_2$.

Now define the action of $W_{w_1 I(w_1, w_2, c, c')}$ on $(w_1 W_{I(w_1, w_2, c, c')}, w_2)$ by $$x \cdot (w_1 v, w_2)=\bigl(x w_1 v \d \i w_2 \d'(x) \i w_2, w_2 \bigr).$$
Then the inclusion map $(w_1 W_{I(w_1, w_2, c, c')}, w_2) \rightarrow [w_1, w_2, c, c']$ induces a bijection between the $W_{w_1 I(w_1, w_2, c, c')}$-orbits on $(w_1 W_{I(w_1, w_2, c, c')}, w_2)$ and the $W_{c'} \times W_c$-cosets in $[w_1, w_2, c, c']$.

Part (2) is proved.

\begin{th7} Each double coset in $W_{c'} \backslash (W_1 \times W_2)/W_c$ contains at most one element of the form $(w_1, w_2)$ with $w_1 \in W^{J_1}$ and $w_2 \in {}^{J'_2} W$.
\end{th7}

Let $(w_1, w_2), (w'_1, w'_2) \in {}^{J'_1} W_1 \times W_2^{J_2}$. By part (1) of the previous Proposition, $W_{c'} (w_1, w_2) W_c \subset [w_1, w_2, c, c']$ and $W_{c'} (w'_1, w'_2) W_c \subset [w'_1, w'_2, c, c']$. In particular, if $W_{c'} (w_1, w_2) W_c=W_{c'} (w'_1, w'_2) W_c$, then $[w_1, w_2, c, c'] \cap [w'_1, w'_2, c, c'] \neq \varnothing$. By the definition, $(w_1, w_2)=(w'_1, w'_2)$. The corollary is proved.

\

It is also worth mentioning the following consequence.

\begin{th8} Let $(W, I)$ be a Coxeter group. Let $J, J' \subset I$ and
$\d: W_J \rightarrow W_{J'}$ be an automorphism with $\d(J)=J'$.
Define the action of $W_J$ on $W$ by $x \cdot y=x y \d(x) \i$. For
$w \in W^{J'}$, set $$I(w, \d)=\max\{K \subset J'; \d w K=K\},
\quad [w, \d]=W_J \cdot (w W_{I(w, \d)}).$$ Then

(1) $W=\sqcup_{w \in W^{J'}} [w, \d]$.

(2) For $w \in W^{J'}$, define an automorphism $\s: W_{I(w, \d)}
\rightarrow W_{I(w, \d)}$ by $\s(v)=\d (w v w \i)$. Then map $W_{I(w,
\d)} \rightarrow W$ defined by $v \rightarrow w v$ induces a
bijection between the $\s$-twisted conjugacy classes in $W_{I(w,
\d)}$ and the $W_J$-orbits in $[w, \d]$.
\end{th8}

Let $(W_1, I_1)=(W_2, I_2)=(W, I)$, $c=(J, J',
\d)$ and $c'=(I, I, id)$. Then the map $W_1 \times W_2 \rightarrow W$ defined by $(w_1,
w_2) \mapsto w_1 \i w_2$ induces a natural bijection $W_{c'}
\backslash (W_1 \times W_2)/W_c$ to the $W_J$-orbits on $W$. Now the
corollary follows easily from Proposition 2.4.

\section{Minimal length elements}

\subsection{} We follow the notation in \cite[section 3.2]{GP2}.

Let $(W, I)$ be a Coxeter group. Let $J, J' \subset
I$ and $\d: W_J \rightarrow W_{J'}$ be an automorphism with
$\d(J)=J'$. Given $w, w' \in W$ and $j \in J$, we write $w
\xrightarrow{s_j}_{\d} w'$ if $w'=s_j w \d(s_j)$ and $l(w') \le
l(w)$. If $w=w_0, w_1, \cdots, w_n=w'$ is a sequence of elements in
$W$ such that for all $k$, we have $w_{k-1} \xrightarrow{s_j}_{\d}
w_k$ for some $j \in J$, then we write $w \rightarrow_{\d} w'$.

We call $w, w' \in W$ {\it elementarily strongly $\d$-conjugate} if $l(w)=l(w')$ and there exists $x \in W_J$ such that $w'=x w \d(x) \i$ and either $l(x w)=l(x)+l(w)$ or $l(w \d(x) \i)=l(x)+l(w)$. We call $w, w'$ {\it strongly $\d$-conjugate} if there is a sequence $w=w_0, w_1, \cdots, w_n=w'$ such that $w_{i-1}$ is elementarily strongly $\d$-conjugate to $w_i$. We will write $w \sim_{\d} w'$ if $w$ and $w'$ are strongly $\d$-conjugate.

If $w \sim_{\d} w'$ and $w \rightarrow_{\d} w'$, then we say that $w$ and $w'$ are in the same $\d$-cyclic shift class and write $w \approx_{\d} w'$. For $w \in W$, set $$\Cyc_{\d}(w)=\{w' \in W; w \approx_{\d} w'\}.$$

If $w' \in \Cyc_{\d}(w)$ for all $w' \in W$ with $w \rightarrow_{\d} w'$, then we call the $\d$-cyclic shift class {\it terminal}. It is easy to see that if $w$ is an element of minimal length in $\{x w \d(x) \i; x \in W_J\}$, then $\Cyc_{\d}(w)$ is terminal.

The following result is proved in \cite{GP1} for the usual conjugacy classes and in \cite{GKP} for the twisted conjugacy classes.

\begin{th9} Let $(W, I)$ be a finite Coxeter group and $\d: W \rightarrow
W$ be an automorphism with $\d(I)=I$. Let $\co$ be a $\d$-twisted
conjugacy class in $W$ and $\co_{\min}$ be the set of minimal length
elements in $\co$. Then

(1) For each $w \in \co$, there exists $w' \in \co_{\min}$ such that
$w \rightarrow_{\d} w'$.

(2) Let $w, w' \in \co_{\min}$, then $w \sim_{\d} w'$.
\end{th9}

\subsection*{3.3} Let $l_1$ (resp. $l_2$) be the length function on
$W_1$ (resp. $W_2$). Define the length function $l$ on $W_1 \times
W_2$ by $l(w_1, w_2)=l_1(w_1)+l_2(w_2)$ for $w_1 \in W_1$ and $w_2
\in W_2$. For each double coset $\co$ in $W_{c'} \backslash (W_1
\times W_2)/W_c$, we set $$\co_{\min}=\{\mathbf w \in \co;
l(\mathbf w) \le l(\mathbf w') \text{ for all } \mathbf w'
\in \co\}.$$

Now following the convention of Fomin and Zelevinsky, we consider
$W_{J'_1} \times W_{J_1}$ as a Coxeter group with simple reflections
$s_{-i}$ (for $i \in J'_1$) and $s_j$ (for $j \in J_1$).

Given $\mathbf w, \mathbf w' \in W_1 \times W_2$ and $i \in
J'_1$, we write $\mathbf w \xrightarrow{s_{-i}}_{c, c'}
\mathbf w'$ if $\mathbf w=(s_i, \d'(s_i)) \mathbf w'$ and
$l(\mathbf w') \le l(\mathbf w)$.

Similarly, given $j \in J_1$, we write $\mathbf w
\xrightarrow{s_j}_{c, c'} \mathbf w'$ if $\mathbf w=\mathbf
w' (s_j, \d(s_j))$ and $l(\mathbf w') \le l(\mathbf w)$.

If $\mathbf w=\mathbf w_0, \mathbf w_1, \cdots, \mathbf
w_n=\mathbf w'$ is a sequence of elements in $W_1 \times W_2$
such that for all $k$, we have $\mathbf w_{k-1}
\xrightarrow{s_i}_{c, c'}  \mathbf w_k$ for some $i \in -J'_1
\sqcup J_1$, then we write $\mathbf w \rightarrow_{c,
c'} \mathbf w'$.

We write $\mathbf w \sim_{c, c'} \mathbf w'$ if there exists a
sequence $\mathbf w=\mathbf w_0, \mathbf w_1, \cdots,
\mathbf w_n=\mathbf w'$ such that $$l(\mathbf
w_{k+1})=l(\mathbf w_k), \quad \mathbf w_{k+1}=(x_k, \d'(x_k))
\mathbf w_k (y_k, \d(y_k))$$ and either $$l \bigl( (x_k, 1)
\mathbf w_k (1, \d(y_k)) \bigr)=l_1(x_k)+l(\mathbf
w_k)+l_1(y_k)$$ or $$l \bigl( (1, \d'(x_k)) \mathbf w_k (y_k, 1)
\bigr)=l_1(x_k)+l(\mathbf w_k)+l_1(y_k)$$ for all $k$ and some
$x_k \in W_{J'_1}$, $y_k \in W_{J_1}$.

We write $\mathbf w \approx_{c, c'} \mathbf w'$ if $\mathbf w \rightarrow_{c, c'} \mathbf w'$ and $\mathbf w \sim_{c, c'} \mathbf w'$.

\begin{th10} Let $(w_1, w_2) \in {}^{J'_1} W_1 \times W_2^{J_2}$ and
$\co$ is a double coset in $W_{c'} \backslash [w_1, w_2, c, c']/W_c$
that corresponds to the $\s$-twisted conjugacy class $\co'$ in
$W_{I(w_1, w_2, c, c')}$ via the map in Proposition 2.4 (2).
Let $\co'_{\min}$ be the set of minimal length elements in $\co'$.
Then

(1) For each $\mathbf w \in \co$, there exists $v \in \co'$ such
that $\mathbf w \rightarrow_{c, c'} (w_1 v, w_2)$.

(2) If $\mathbf w \in \co_{\min}$, then there exists $v \in
\co'_{\min}$ such that $\mathbf w \approx_{c, c'} (w_1 v, w_2)$.
\end{th10}

Let $(J_1^{(n)}, J_2^{' (n)}, w_1^{(n)}, w_2^{(n)}, u_1^{(n)},
u_2^{(n)}, v_1^{(n)}, v_2^{(n)})_{n \ge 0}$ be the sequence
associated to $\mathbf w=(w'_1, w'_2)$. Then it is easy to see
that
$$\mathbf w \rightarrow_{c, c'} (u_1^{(0)}, \d'(v_1^{(0)}) \i w'_2)
\rightarrow_{c, c'} (u_1^{(0)} \d \i(v_2^{(0)}) \i, u_2^{(0)})$$ and
for $n \ge 0$, \begin{align*} (u_1^{(n)} \d \i(v_2^{(n)}) \i,
u_2^{(n)}) & \rightarrow_{c, c'} (u_1^{(n+1)}, \d' (v_1^{(n+1)}) \i
u_2^{(n)}) \\ & \rightarrow_{c, c'} (u_1^{(n+1)} \d \i(v_2^{(n+1)}) \i,
u_2^{(n+1)}). \end{align*}

By the proof of Proposition 2.4, $u_1^{(n)}=w_1$,
$u_2^{(n)}=w_2$ and $\d \i (v_2^{(n)}) \in \co'$ for $n \gg 0$. Thus
$\mathbf w \rightarrow_{c, c'} (w_1 v, w_2)$ for some $v \in
\co'$.

Part (1) is proved.

If moreover, $\mathbf w \in \co_{\min}$, then $(w_1 v, w_2) \in
\co_{\min}$. By Proposition 2.4 (2), $v \in \co'_{\min}$. It is
then easy to see that $$\mathbf w \approx_{c, c'} (u_1^{(0)},
\d'(v_1^{(0)}) \i w'_2) \approx_{c, c'} (u_1^{(0)} \d \i(v_2^{(0)}) \i,
u_2^{(0)})$$ and for $n \ge 0$, \begin{align*} (u_1^{(n)} \d
\i(v_2^{(n)}) \i, u_2^{(n)}) & \approx_{c, c'} (u_1^{(n+1)}, \d'
(v_1^{(n+1)}) \i u_2^{(n)})
\\ & \approx_{c, c'} (u_1^{(n+1)} \d \i(v_2^{(n+1)}) \i, u_2^{(n+1)}).
\end{align*}

In particular, $\mathbf w \approx_{c, c'} (w_1 v, w_2)$. Part (2) is
proved.

\

Now combining the above Proposition with Theorem 3.2, we have
the following consequence.

\begin{th11} Let $(w_1, w_2) \in {}^{J'_1} W_1 \times W_2^{J_2}$ with
$W_{I(w_1, w_2, c, c')}$ is a finite Coxeter group. Let $\co \in W_{c'} \backslash
[w_1, w_2, c, c']/W_c$. Then

(1) For each $\mathbf w \in \co$, there exists $\mathbf w' \in
\co_{\min}$ such that $\mathbf w \rightarrow_{c, c'} \mathbf
w'$.

(2) Let $\mathbf w, \mathbf w' \in \co_{\min}$, then
$\mathbf w \sim_{c, c'} \mathbf w'$.
\end{th11}

\begin{rmk} By definition, if $W_{J_1}$ or $W_{J'_1}$ is a finite
Coxeter group, then $W_{I(w_1, w_2, c, c')}$ is also a finite
Coxeter group.
\end{rmk}

\begin{th12} Let $(W, I)$ be a Coxeter group and $\d: W \rightarrow
W$ be an automorphism with $\d(I)=I$. Let $w \in W$ with $\d(w)=w \i$. Then $w \rightarrow_{\d} w_J$ for some $J=\d(J) \subset I$. Moreover, $w_J \d(s_j)=s_j w_J$ for $j \in J$ and $w_J$ has minimal length in its $\s$-conjugacy class in $W$. In particular, $x \d(x) \i
\rightarrow_{\d} 1$ for all $x \in W$.
\end{th12}

\begin{rmk} This is a generalization of Richardson's theorem in \cite{R}. Our proof is similar to the proof of \cite[3.2.10]{GP2} which was essentially due to Howlett.
\end{rmk}

We argue by induction on $l(w)$. For $w=1$ this is clear. Suppose that $l(w) \ge 1$. Since $\d(w)=w \i$, we have that $\{i \in I; l(s_i w)<l(w)\}=\{i \in I; l(w
\d(s_i))<l(w)\}$. Set $$J=\{i \in I; s_i w=w \d(s_i), l(s_i
w)<l(w)\}.$$ Then $w=w_J w'$, where $w_J$ is the maximal element in
$W_K$ and $w' \in {}^J W$.

If $w'=1$, then $w=w_J=\d(w_J) \i$ and $\d(J)=J$. Now for $x=a
b$ with $a \in W^J$ and $b \in W_J$, $x w \d(x) \i=a \bigl( b w \d(b) \i
\bigr) \d(a) \i=a w_J \d(a) \i$. So $l(x w \d(x) \i) \ge l(a w_J)-l(a)=l(w)$. 

If $w' \neq 1$, then there exists $i \in I$ with $l(w' \d(s_i))<l(w')$. Hence
$l(w \d(s_i))<l(w)$ and $l(s_i w)<l(w)$. If $i \in J$, then $s_i w
\in W_J w'$, $w \d(s_i) \in W_J w' \d(s_i)$ and $w', w' \d(s_i) \in
{}^J W$. Hence $s_i w \neq w \d(s_i)$. That is a contradiction.
Hence $i \notin J$. By \cite[Lemma 1.2.6]{GP2}, $l(s_i w
\d(s_i))<l(w)$. Hence $w \rightarrow_{\d} s_i w \d(s_i)$. Now the proposition follows from 
induction hypothesis.

\

Now combining the above Proposition with Proposition 3.4, we have
the following consequence.

\begin{th13} Let $(w_1, w_2) \in {}^{J'_1} W_1 \times W_2^{J_2}$ and
$\co=W_{c'} (w_1, w_2) W_c$. Let $\mathbf w \in \co$. Then
$\mathbf w \rightarrow_{c, c'} (w_1, w_2)$. If moreover,
$\mathbf w \in \co_{\min}$, then $\mathbf w \approx_{c, c'} (w_1,
w_2)$.
\end{th13}

\

It is also worth mentioning the following consequence which is a
generalization of Theorem 3.2.

\begin{th14} Let $(W, I)$ be a Coxeter group. Let $J, J' \subset I$ and
$\d: W_J \rightarrow W_{J'}$ be an automorphism with $\d(J)=J'$.
Define the action of $W_J$ on $W$ by $x \cdot y=x y \d(x) \i$. Let
$\co$ be a $W_J$-orbit in $W$ and $\co_{\min}$ be the set of minimal
length elements in $\co$. If moreover, $W_J$ is a finite Coxeter
group or $\co \cap W^{J'} \neq \varnothing$. Then

(1) For each $w \in \co$, there exists $w' \in \co_{\min}$ such that
$w \rightarrow_{\d} w'$.

(2) Let $w, w' \in \co_{\min}$, then $w \sim_{\d} w'$.
\end{th14}

\begin{rmk} The case when $W$ is a finite Coxeter group was proved in \cite{F}.
\end{rmk}

\section{Distinguished double cosets}

\begin{th15} For $(w_1, w_2) \in W_1 \times W_2$, set \begin{align*} &
W(w_1, w_2)=\{v \in W_{J_1}; \d' \Ad(w_1) \bigl(\d \i \Ad(w_2) \i
\d' \Ad(w_1) \bigr)^n v \in W_{J'_2} \\ & \text{ and } \bigl(\d \i
\Ad(w_2) \i \d' \Ad(w_1) \bigr)^{n+1} v \in W_{J_1} \text{ for all }
n \ge 0\}.\end{align*}

Then \begin{gather*} \tag{1} W(w_1, w_2)=W_{I(w_1, w_2, c, c')}
\text{ for } (w_1, w_2) \in {}^{J'_1} W_1 \times W_2^{J_2}. \\
\tag{2} W(w'_1, w'_2)=W_{\d \i I(w'_2, w'_1, c \i, (c') \i)} \text{
for } (w'_1, w'_2) \in W_1^{J_1} \times {}^{J'_2} W_2. \end{gather*}
\end{th15}

Part (2) is equivalent to part (1). So we will only prove part (1).

Let $(J_1^{(n)}, J_2^{(n)}, w_1^{(n)}, w_2^{(n)})_{n \ge 0}$ be the
element in $\ct(c, c')$ whose image under the map $\phi$ defined in
Proposition 1.7 is $(w_1, w_2)$. Let $v \in W(w_1, w_2)$. Then
we can prove by induction on $n \ge 0$ that
\begin{gather*} \d' \Ad(w_1) \bigl(\d
\i \Ad(w_2) \i \d' \Ad(w_1) \bigr)^n v \in W_{J_2^{(n)}}, \\
\bigl(\d \i \Ad(w_2) \i \d' \Ad(w_1) \bigr)^{n+1} v \in
W_{J_1^{(n+1)}}.
\end{gather*}

In particular, for $n \gg 0$, $$\bigl(\d \i \Ad(w_2) \i \d' \Ad(w_1)
\bigr)^{n+1} v \in W_{I(w_1, w_2, c, c')}.$$ Hence $v \in
W_{I(w_1, w_2, c, c')}$. On the other hand $W_{I(w_1, w_2, c, c')}
\subset W(w_1, w_2)$. Hence $W_{I(w_1, w_2, c, c')}=W(w_1, w_2)$.

Part (1) is proved.

\subsection*{4.2} A double coset $\co$ in $W_{c'} \backslash (W_1 \times W_2)/W_c$
is called {\it distinguished} with respect to $c, c'$ if it contains
some element $(w_1, w_2)$ with $w_1 \in {}^{J'_1} W_1$ and $w_2 \in
W_2^{J_2}$. In this case, we simply write $I(\co, c, c')$ for
$I(w_1, w_2, c, c')$ and $[\co, c, c']$ for $[w_1, w_2, c, c']$.

The minimal length elements in distinguished double cosets are
called {\it distinguished elements} in $W_1 \times W_2$ with respect to
$c, c'$.

\begin{th16} Each element in $W_1^{J_1} \times {}^{J'_2} W_2$ is a
distinguished element in $W_1 \times W_2$ with respect to $c, c'$.
\end{th16}

Let $w'_1 \in W_1^{J_1}$ and $w'_2 \in {}^{J'_2} W_2$. Then $(w'_1,
w'_2)$ is of minimal length in $W_{c'} (w'_1, w'_2) W_c$. By
Proposition 2.4, we may assume that $(w'_1, w'_2)=(x, \d'(x))
(w_1 v_1, w_2 v_2) (y, \d(y))$ for some $x \in W_{J'_1}$, $y \in W_{J_1}$,
$w_1 \in {}^{J'_1} W_1$, $w_2 \in W_2^{J_2}$, $v_1 \in W_{I(w_1,
w_2, c, c')}$ and $v_2 \in W_{\d I(w_1,
w_2, c, c')}$. It is easy to see that we may assume furthermore that
$x \in W_{J'_1}^{I'(w_1, w_2, c, c')}$ and $y \in W_{J_1}^{I(w_1,
w_2, c, c')}$, where $I'(w_1, w_2, c, c')=\d' w_1 I(w_1, w_2, c,
c')$.

By Lemma 4.1, \begin{align*} W_{\d \i I(w'_2, w'_1, c \i, (c') \i)} &=W(w'_1, w'_2)=y \i W(w_1, w_2) y \\ &=y \i
W_{I(w_1, w_2, c, c')} y.\end{align*}

Hence $\d \i I(w'_2, w'_1, c \i, (c') \i)=y
\i I(w_1, w_2, c, c')$. Thus $x v_1 y \in W_{\d \i I(w'_2, w'_1, c \i, (c') \i)}$ and
$x w_1 y (y \i v y)=w'_1 \in W_1^{\d \i I(w'_2, w'_1, c \i, (c')
\i)}$. Hence $y \i v_1 y=1$ and $v_1=1$. Similarly, $v_2=1$. Then $(w'_1, w'_2) \in W_{c'} (w_1, w_2) W_c$ and $(w'_1, w'_2)$ is a distinguished element.

The proposition is proved.

\

In the rest of this section, we will introduce a partial order on the set of distinguished double cosets.

\begin{th17} Let $\co \in W_{c'} \backslash (W_1 \times W_2)/W_c$ and $\mathbf w \in \co_{\min}$. If $\mathbf w \le \mathbf w'$ and $\mathbf w'_1 \rightarrow_{c, c'} \mathbf w'$, then there exists $\mathbf w_1 \in \co_{\min}$ with $\mathbf w_1 \le \mathbf w'_1$.
\end{th17}

It suffices to prove that for any $i \in -J'_1 \sqcup J_1$ and $\mathbf w'_1 \in W_1 \times W_2$ with $\mathbf w'_1 \xrightarrow{s_i}_{c, c'} \mathbf w'$, there exists $\mathbf w_1 \in \co_{\min}$ with $\mathbf w_1 \le \mathbf w'_1$.

There are two cases:

Case 1: $\mathbf w'_1=(s_i, \d'(s_i)) \mathbf w'$ for $i \in J'_1$;

Case 2: $\mathbf w'_1=\mathbf w' (s_i, \d(s_i))$ for $i \in J_1$.

We will only prove for Case 1. Case 2 can be proved in the same way.

Assume that $\mathbf w'=(w'_1, w'_2)$ and $\mathbf w=(w_1, w_2)$. If $l(\mathbf w_1)>l(\mathbf w')$. Then $\mathbf w \le \mathbf w' \le \mathbf w'_1$. Now assume that $l(\mathbf w_1)=l(\mathbf w')$. Without loss of generalization, we may assume that $s_i w'_1 <w'_1$ and $\d'(s_i) w'_2>w'_2$.

If $s_i w_1>w_1$, then by \cite[Corollary 2.5]{L1}, we have that $w_1 \le s_i w'_1$. We also have that $w_2 \le w'_2<\d'(s_i) w'_2$. Hence $\mathbf w \le \mathbf w'_1$.

If $s_i w_1<w_1$, then by \cite[Corollary 2.5]{L1}, we have that $s_i w_1<s_i w'_1$. Since $\d'(s_i) w'_2>w'_2$ and $w_2 \le w'_2$, then by \cite[Corollary 2.5]{L1}, we also have that $\d'(s_i) w_2 \le \d'(s_i) w'_2$. Hence $(s_i, \d'(s_i)) \mathbf w \le \mathbf w'_1$. Moreover, since $s_i w_1<w_1$, $l((s_i, \d'(s_i)) \mathbf w) \le l(\mathbf w)$ and  $\mathbf w \in \co_{\min}$, we have that $(s_i, \d'(s_i)) \mathbf w \in \co_{\min}$.

The lemma is proved.

\

Now combining the above lemma with Corollary 3.8, we have the following consequence.

\begin{th18} Let $\co \in W_{c'} \backslash (W_1 \times W_2)/W_c$. If moreover, any of the following condition holds:

(1) $W_{J'}$ or $W_{J_1}$ is a finite Coxeter group;

or (2) $\co$ is a distinguished double coset,

then $\co_{\min}=\{\mathbf w \in \co; \mathbf w \text{ is a minimal element in } \co\}$.
\end{th18}

\

Here is another consequence.

\begin{th19} Let $\co, \co'$ be distinguished double cosets in $W_{c'} \backslash (W_1 \times W_2)/W_c$. Then the following conditions are equivalent:

(1) For some $\mathbf w' \in \co'_{\min}$, there exists $\mathbf w \in \co_{\min}$ such that $\mathbf w \le \mathbf w'$.

(2) For any $\mathbf w' \in \co'_{\min}$, there exists $\mathbf w \in \co_{\min}$ such that $\mathbf w \le \mathbf w'$.
\end{th19}

\subsection*{4.7} Now we define a partial order on the set of distinguished double cosets in $W_{c'} \backslash (W_1 \times W_2)/W_c$ as follows:

$\co \le \co'$ if for some (or equivalently, any) $\mathbf w' \in \co'_{\min}$, there exists $\mathbf w \in \co_{\min}$ with $\mathbf w \le \mathbf w'$.

This partial order will be used in the next section to describe the closure relations of the so-called $\car_{\cc'} \times \car_{\cc}$-stable pieces.

\section{$\car_{\cc'} \times \car_{\cc}$-stable pieces and Hecke algebras}

\subsection{} For $i=1, 2$, let $G_i$ be a connected reductive algebraic group over an algebraically closed field $k$, $B_i$ be a Borel subgroup of $G_i$ and $T_i \subset B_i$ be a maximal torus. $B_i$ and $T_i$ determine a Weyl group $W_i$ and the set $I_i$ of its simple reflections. For $w \in W$, we use the same symbol $w$ for a representative of $w$ in $N(T)$. For each subset $J_i$ of $I_i$, we denote by $P_{J_i}$ the standard parabolic subgroup of type $J_i$, $L_{J_i}$ the Levi subgroup of $P_{J_i}$
that contains $T_i$ and $\pi_{J_i}: P_{J_i} \rightarrow L_{J_i}$ the projection map.

For any subvariety $X$ of $G_1 \times G_2$, we denote by $\bar{X}$ its closure in $G_1 \times G_2$.

{\it An admissible triple} of $G_1 \times G_2$ is by definition a triple $\cc=(J_1, J_2, \th_{\d})$ consisting of $J_1 \subset I_1$, $J_2 \subset I_2$, an isomorphism $\d: W_{J_1} \rightarrow W_{J_2}$ with $\d(J_1)=J_2$ and an isomorphism $\th_{\d}: L_{J_1} \rightarrow L_{J_2}$ that maps $T_1 \subset L_{J_1}$ to $T_2 \subset L_{J_2}$ and the root subgroup $U_{\a_i} \subset L_{J_1}$ to the root subgroup $U_{\a_{\d(i)}} \subset L_{J_2}$ for $i \in J_1$. To each admissible triple $\cc=(J_1, J_2, \th_{\d})$,
we associate a subgroup $\car_{\cc}$ of $G_1 \times G_2$ defined as follows $$\car_{\cc}=\{(p, q); p \in P_{J_1}, q \in P_{J_2}, \th_{\d}(\pi_{J_1}(p))=\pi_{J_2}(q)\}.$$ Moreover, each admissible triple $\cc=(J_1, J_2, \th_{\d})$ of $G_1 \times G_2$ determines an admissible triple $c=(J_1, J_2, \d)$ of $W_1 \times W_2$. We also set $\cb_{\cc}=\car_{\cc} \cap (B_1, B_2)$.

Notice that if $G_1=G_2=G$, $B_1=B_2=B$, $T_1=T_2=T$ and $I_1=I_2=I$, then $\car_{(I, I, id)}$ is the diagonal subgroup $G_{\D}$ of $G \times G$ and $\cb_{\cc}=B_{\D}$.

In fact, in \cite{LY}, they consider a slightly more general class of the groups $\car_{\cc}$. However, the results below can be easily generalized to the more general setting.

\subsection{} Now given admissible triples $\cc=(J_1, J_2, \th_{\d})$ and $\cc'=(J'_1, J'_2, \th_{\d'})$ of $G_1 \times G_2$. Let $c=(J_1, J_2, \d)$ and $c'=(J'_1, J'_2, \d')$ be the corresponding admissible triples of $W_1 \times W_2$. For $(w_1, w_2) \in {}^{J'_1} W_1 \times W_2^{J_2}$, define $$[w_1, w_2, \cc, \cc']=\car_{\cc'} (B_1 w_1 B_1, B_2 w_2 B_2) \car_{\cc}.$$

We call $[w_1, w_2, \cc, \cc']$ a $\car_{\cc'} \times \car_{\cc}$-stable piece of $G_1 \times G_2$.

\begin{th20} Let $(w'_1, w'_2)$ be a distinguished element in $W_1 \times W_2$ with respect to $c, c'$ and $\pi(w'_1, w'_2)=(w_1, w_2)$. Then $$[w_1, w_2, \cc, \cc']=\car_{\cc'} (B_1 w'_1 B_1, B_2 w'_2 B_2) \car_{\cc'}.$$
\end{th20}

\begin{rmk} Thus for a distinguished double $\co$ in $W_{c'} \backslash (W_1 \times W_2)/W_c$, we may write $[\co, \cc, \cc']$ for $[w_1, w_2, \cc, \cc']$ where $(w_1, w_2)$ is the unique element in $\co \cap ({}^{J'_1} W_1 \times W_2^{J_2})$.
\end{rmk}

Let $(J_1^{(n)}, J_2^{' (n)}, w_1^{(n)}, w_2^{(n)}, u_1^{(n)}, u_2^{(n)}, v_1^{(n)}, v_2^{(n)})_{n \ge 0}$ be the sequence associated to $(w'_1, w'_2)$. By the proof of Proposition 3.4, \begin{align*} \car_{\cc'} (B_1 w'_1 B_1, B_2 w'_2 B_2) \car_{\cc'} &=\car_{\cc'} (B_1 u_1^{(0)} B_1, B_2 \d'(v_1^{(0)}) \i w'_2 B_2) \car_{\cc'} \\ &=\car_{\cc'} (B_1 u_1^{(0)} \d \i(v_2^{(0)}) \i B_1, B_2 u_2^{(0)} B_2) \car_{\cc'}\end{align*} and for $n \ge 0$, \begin{align*} & \car_{\cc'} (B_1 (u_1^{(n)} \d
\i(v_2^{(n)}) \i B_1, B_2 u_2^{(n)} B_2) \car_{\cc'} \\ &=\car_{\cc'} (B_1 u_1^{(n+1)} B_1, B_2 \d' (v_1^{(n+1)}) \i u_2^{(n)} B_2) \car_{\cc'} \\ &=\car_{\cc'} (B_1 u_1^{(n+1)} \d \i(v_2^{(n+1)}) \i B_1, B_2 u_2^{(n+1)} B_2) \car_{\cc'}.\end{align*}

By the proof of Proposition 2.4, $u_1^{(n)}=w_1$,
$u_2^{(n)}=w_2$ for $n \gg 0$. Moreover, since $l(w'_1, w'_2)=l(w_1, w_2)$, we have that $v_2^{(n)}=1$ for $n \gg 0$. Thus
$[w_1, w_2, \cc, \cc']=\car_{\cc'} (B_1 w'_1 B_1, B_2 w'_2 B_2) \car_{\cc'}$.

The proposition is proved.

\

Now combining the above proposition with Proposition 4.3, we have the following consequence.

\begin{th21} $\car_{\cc'} (B_1 w_1 B_1, B_2 w_2 B_2) \car_{\cc}$ is a $\car_{\cc'} \times \car_{\cc}$-stable piece for $(w_1, w_2) \in W_1^{J_1} \times {}^{J'_2} W_2$.
\end{th21}

\begin{rmk} In \cite{LY}, Lu and Yakimov define the $\car_{\cc'} \times \car_{\cc}$-stable piece using $W_1^{J_1} \times {}^{J'_2} W_2$ instead of $^{J'_1} W_1 \times W_2^{J_2}$. Now we can see from the above corollary that our definition coincide with theirs.
\end{rmk}

\

We may reformulate the above corollary in a different way.

\begin{th22} Let $\partial: G_1 \times G_2 \rightarrow G_2 \times G_1$ be the map defined by $(g_1, g_2) \mapsto (g_2, g_1)$. Then $\partial$ sends a $\car_{\cc'} \times \car_{\cc}$-stable piece of $G_1 \times G_2$ to a $\car_{(\cc') \i} \times \car_{\cc \i}$-stable piece of $G_2 \times G_1$.
\end{th22}

\begin{rmk}
A special case of the corollary has been proved in \cite[Proposition 2.5]{H3}. The proof here is simpler.
\end{rmk}

\

We also have the following properties of the $\car_{\cc'} \times \car_{\cc}$-stable pieces which were proved in \cite{LY} generalizing some results of the $G$-stable pieces obtained in \cite{L3}.

\begin{th23} (1) $G_1 \times G_2=\sqcup_{(w_1, w_2) \in ^{J'_1} W_1 \times W_2^{J_2}} [w_1, w_2, \cc, \cc']$.

(2) Let $(w_1, w_2) \in ^{J'_1} W_1 \times W_2^{J_2}$. Define an automorphism $\th_{\s}: L_{I(w_1, w_2, c, c')} \rightarrow
L_{I(w_1, w_2, c, c')}$ by $\th_{\s}(l)=\th_{\d} \i \bigl(w_2 \i \th_{\d} (w_1 l w_1 \i) w_2 \bigr)$. Then map $L_{I(w_1, w_2, c, c')} \rightarrow G_1
\times G_2$ defined by $l \rightarrow (w_1 l, w_2)$ induces a
bijection between the $\th_{\s}$-twisted conjugacy classes on $L_{I(w_1,
w_2, c, c')}$ and the double cosets $\car_{\cc'} \backslash [w_1, w_2, \cc, \cc'] /\car_{\cc}$.
\end{th23}

\begin{rmk} This proposition is an analogy of Proposition 2.4. In fact, we can prove this proposition using a modified version of the inductive method in 2.1. The case for the $G$-stable pieces was showed in this way in \cite[4.3 \& 4.4]{H1}.
\end{rmk}

\

The following property of the $\car_{\cc'} \times \car_{\cc}$-stable piece will be used to study the closure relations.

\begin{th24} Let $\mathbf w$ be a distinguished element in $W_1 \times W_2$ with respect to $c, c'$. Then $\car_{\cc'} (T_1, T_2) \mathbf w \car_{\cc}$ is dense in the $\car_{\cc'} \times \car_{\cc}$-stable piece $\car_{\cc'} (B_1, B_2) \mathbf w (B_1, B_2) \car_{\cc}$.
\end{th24}

By Corollary 3,7, it suffices to prove the case where $\mathbf w=(w_1, w_2) \in {}^{J'_1} W_1 \times W_2^{J_2}$. In this case, by part (2) of the previous Proposition, $$\car_{\cc'} (B_1, B_2) \mathbf w (B_1, B_2) \car_{\cc}=\car_{\cc'} \mathbf w (L_{I(w_1, w_2, c, c')}, 1) \car_{\cc}.$$

Let $\th_{\s}: L_{I(w_1, w_2, c, c')} \rightarrow
L_{I(w_1, w_2, c, c')}$ be the automorphism defined in part (2) of the previous Proposition. Then $$\car_{\cc'} (T_1, T_2) \mathbf w \car_{\cc}=\car_{\cc'} \mathbf w (L', 1) \car_{\cc},$$ where $L'=\{l t \th_{\s} (l) \i; l \in L_{I(w_1, w_2, c, c')}, t \in T\}$. By \cite[Lemma 4]{Sp}, $L'$ is dense in $L_{I(w_1, w_2, c, c')}$. Hence $\car_{\cc'} (T_1, T_2) \mathbf w \car_{\cc}$ is dense in $[w_1, w_2, c, c']$. The lemma is proved.

\begin{th25} Let $\mathbf w \in W_1 \times W_2$, then $$\overline{\car_{\cc'} (B_1, B_2) \mathbf w (B_1, B_2) \car_{\cc}}=\sqcup_{\co} [\co, \cc, \cc']$$ where $\co$ runs over the distinguished double cosets in $W_{c'} \backslash (W_1 \times W_2)/W_c$ that contains a minimal length element $\mathbf w'$ with $\mathbf w' \le \mathbf w$.
\end{th25}

\begin{rmk} This was first proved in \cite[Theorem 5.2]{LY}, which is a generalization of \cite[Corollary 5.5]{H2}.
\end{rmk}

We will simply write $\car$ for $\car_{\cc'} \times \car_{\cc}$ and $\cb$ for $\cb_{\cc'} \times \cb_{\cc}$. Define the action of $\cb_{\cc'} \times \cb_{\cc}$ on $\car_{\cc'} \times \car_{\cc} \times (G_1 \times G_2)$ by $(\mathbf b_1, \mathbf b_2) \cdot (\mathbf g_1, \mathbf g_2, \mathbf g)=(\mathbf g_1 \mathbf b_1 \i, \mathbf g_2 \mathbf b_2 \i, \mathbf b_1 \mathbf g \mathbf b_2 \i)$. Let $\car \times_{\cb} (G_1 \times G_2)$ be its quotient space. Define the map $\car_{\cc'} \times \car_{\cc} \times (G_1 \times G_2) \rightarrow G_1 \times G_2$ by $(\mathbf g_1, \mathbf g_2, \mathbf g) \mapsto \mathbf
g_1 \mathbf g \mathbf g_2 \i$. Then this map induces a proper map $\car \times_{\cb} (G_1 \times G_2) \rightarrow G_1 \times G_2$. In particular, $$\car_{\cc'} \overline{(B_1, B_2) \mathbf w (B_1, B_2)} \car_{\cc}=\overline{\car_{\cc'} (B_1, B_2) \mathbf w (B_1, B_2) \car_{\cc}}.$$

Now let $\co$ be a distinguished double cosets in $W_{c'} \backslash (W_1 \times W_2)/W_c$ that contains a minimal length element $\mathbf w'$ with $\mathbf w' \le \mathbf w$. Then $$\car_{\cc'} (T_1, T_2) \mathbf w' \car_{\cc} \subset \car_{\cc'} (B_1, B_2) \mathbf w' \car_{\cc} \subset \car_{\cc'} \overline{(B_1, B_2) \mathbf w (B_1, B_2)} \car_{cc}.$$ By the previous lemma, $\car_{\cc'} (T_1, T_2) \mathbf w' \car_{\cc}$ is dense in $[\co, \cc, \cc']$. Thus $[\co, \cc, \cc'] \subset
\overline{\car_{\cc'} (B_1, B_2) \mathbf w (B_1, B_2) \car_{cc}}$.

Now it suffices to prove that $\car_{\cc'} (B_1, B_2) \mathbf w (B_1, B_2) \car_{\cc} \subset \sqcup_{\co} [\co, \cc, \cc']$ where $\co$ runs over the distinguished double cosets in $W_{c'} \backslash (W_1 \times W_2)/W_c$ that contains a minimal length element $\mathbf w'$ with $\mathbf w' \le \mathbf w$.

We argue by induction on $l(\mathbf w)$. For $\mathbf w=1$, the statement is clear. Assume that $l(\mathbf w)>1$. Let $(J_1^{(n)}, J_2^{' (n)}, w_1^{(n)}, w_2^{(n)}, u_1^{(n)}, u_2^{(n)}, v_1^{(n)}, v_2^{(n)})_{n \ge 0}$ be the sequence associated to $\mathbf w$. Then we can prove by induction on $n$ that

\begin{align*} \car_{\cc'} (B_1, B_2) & \mathbf w (B_1, B_2) \car_{\cc} \subset \cup_{\mathbf w'<\mathbf w} \car_{\cc'} (B_1, B_2) \mathbf w' (B_1, B_2) \car_{\cc} \\ & \cup \car_{\cc'} (B_1, B_2) (u_1^{(n)} \d \i (v_2^{(n)}) \i, u_2^{(n)}) (B_1, B_2) \car_{\cc}
\\ & \subset \cup_{\mathbf w' <\mathbf w} \car_{\cc'} (B_1, B_2) \mathbf w' (B_1, B_2) \car_{\cc} \\ & \cup \car_{\cc'} (B_1, B_2) (u_1^{(n+1)}, \d' (v_1^{(n+1)}) \i u_2^{(n)}) (B_1, B_2) \car_{\cc}.
\end{align*}

By induction hypothesis and Proposition 1.7, the statement holds for $\mathbf w$. The Proposition is proved.

\begin{th26} Let $\co$ be a distinguished double coset in $W_{c'} \backslash (W_1 \times W_2)/W_c$. Then $$\overline{[\co, \cc, \cc']}=\sqcup_{\co' \text{ is a distinguished double coset in } W_{c'} \backslash (W_1 \times W_2)/W_c, \co' \le \co} [\co', \cc, \cc'].$$
\end{th26}

\section{Unipotent character sheaves}

\subsection{} We follow the notation of \cite{BBD}. Let $X$ be an
algebraic variety over $\mathbf k$ and $l$ be a fixed prime number
invertible in $\mathbf k$. We write $\cd(X)$ instead of $\cd^b_c(X,
\bar{\mathbb Q}_l)$. If $C \in \cd(X)$ and $A$ is a simple
perverse sheaf on $X$, we write $A \dashv C$ if $A$ is a
composition factor of ${}^p H^i(C)$ for some $i \in \mathbb Z$.
For $A, B \in \cd(X)$, we write $A=B[\cdot]$ if $A=B[m]$ for some
$m \in \mathbb Z$.

Let $C, C_1, \cdots C_n \in \cd(X)$. We write $C \in <C_i; i=1, 2, \cdots, n>$ if there exist $m>n$ and
$C_{n+1}, \cdots, C_m \in \cd(X)$ such that $C_m=C$ and for each
$n+1 \le i \le m$, there exists $1 \le j, k<i$ such that $(C_j[\cdot], C_i,
C_k[\cdot])$ is a distinguished triangle in $\cd(X)$. In this case, if $A \dashv C$, then $A \dashv C_i$ for some $1 \le i \le n$.

Let $H$ be a connected algebraic group and $X$, $Y$ be varieties
with a free $H$-action on $X \times Y$. Denote by $X \times^H Y$
the quotient space. For $C_1 \in \cd(X)$ and $C_2 \in \cd(Y)$ such
that $C_1 \boxtimes C_2$ is $H$-equivariant, we denote by $C_1
\odot C_2$ be the element in $\cd(X \times^H Y)$ whose inverse
image under $X \times Y \rightarrow X \times^H Y$ is $C_1
\boxtimes C_2$.

\subsection{} We keep the notation in 5.1. For $w_1 \in W_1$, we denote by $\cl_{w_1}$ the trivial local system on $B_1 w_1 B_1$. We also use the same notation for its extension by $0$ to $G_1$. Let $\ca_{w_1}$ be its perverse extension to $G_1$, i.e., a perverse sheaf on $G_1$ supported by $\overline{B_1 w_1 B_1}$ and the restriction to $B_1 w_1 B_1$ is $\cl_{w_1}[\dim(B_1 w_1 B_1)]$. We can define $\cl_{w_2}$ and $\ca_{w_2}$ for $w_2 \in W_2$ in the same way. For $\mathbf w=(w_1, w_2) \in W_1 \times W_2$, set $\cl_{\mathbf w}=\cl_{w_1} \boxtimes \cl_{w_2}$ and $\ca_{\mathbf w}=\ca_{w_1} \boxtimes \ca_{w_2}$.

We will simply write $\car$ for $\car_{\cc'} \times \car_{\cc}$ and $\cb$ for $\cb_{\cc'} \times \cb_{\cc}$. Then we have a proper map $\pi: \car \times_{\cb} (G_1 \times G_2) \rightarrow G_1 \times G_2$. See the proof of Proposition 5.8.

We call a simple perverse sheaf $C$ on $G_1 \times G_2$ a {\it unipotent character sheaf} with respect to $\cc$ and $\cc'$ if $C$ is a constitute of $\pi_!(\bar{\mathbb Q}_l[\dim(\car)] \odot \ca_{\mathbf w})$ for some $\mathbf w \in W_1 \times W_2$. This is a generalization of Lusztig's unipotent parabolic character sheaves in \cite{L3}.

We may also define character sheaves with respect to $\cc$ and $\cc'$ by using tame local systems instead of trivial local systems. However, we will not go into details here.

Recently, Springer told me that he also got a similar generalization of Lusztig's parabolic character sheaves.

\begin{th27} Let $\mathbf w, \mathbf w' \in W_1 \times W_2$. Then 

(1) If $\mathbf w \rightarrow_{c, c'} \mathbf w'$ and $l(\mathbf w)>l(\mathbf w')$, then $$\pi_! (\bar{\mathbb Q}_l \odot \cl_{\mathbf w}) \in <\pi_! (\bar{\mathbb Q}_l \odot \cl_{\mathbf x})>_{l(\mathbf x)<l(\mathbf w)}.$$

(2) If $\mathbf w \approx_{c, c'} \mathbf w'$, then $$\pi_! (\bar{\mathbb Q}_l \odot \cl_{\mathbf w})=\pi_! (\bar{\mathbb Q}_l \odot \cl_{\mathbf w'}).$$
\end{th27}

\begin{rmk} The proof is similar to \cite[Lemma 3.9]{H2}.
\end{rmk}

It suffices to prove the case where $\mathbf w_1 \xrightarrow{s_i}_{c, c'} \mathbf w_2$ for some $i \in -J'_1 \sqcup J_1$. Without loss of generalization, we may assume that $i \in J_1$.
We assume that $\mathbf w=(w_1, w_2)$. Then $\mathbf w'=(w_1 s_i, w_2 s_{\d(i)})$. Since $l(\mathbf w) \ge l(\mathbf w')$, either $w_1>w_1 s_i$ or $w_2>w_2 s_{\d(i)}$. We assume that $w_1>w_1 s_i$. The other case can be proved in the same way.

Set $B_{J_i}=B_i \cap L_{J_i}$. For $w \in W_{J_i}$, let $\cl'_{w}$ be the trivial local system on $B_{J_i} w B_{J_i}$. Define the action of $B_{J_i}$ on $G_i \times L_{J_i}$ by $b \cdot (g, g')=(g b \i, b g')$. Let $G_i \times^{B_{J_i}} L_{J_i}$ be the quotient space.

Define the action of $\cb$ on $\car \times \bigl( (G_1 \times^{B_{J_1}} L_{J_1}) \times G_2 \bigr)$ by $b \cdot (r, (g, g'), g_2)=(r b \i, (b \cdot g, g'), b \cdot g_2)$. Let $\car \times_{\cb} \bigl( (G_1 \times^{B_{J_1}} L_{J_1}) \times G_2 \bigr)$ be the quotient. The map $\car \times \bigl( (G_1 \times^{B_{J_1}} L_{J_1}) \times G_2 \bigr) \rightarrow G_1 \times G_2$ defined by $(r, (g, g'), g_2) \mapsto r \cdot (g g', g_2)$ induces a proper morphism $$f_{1, 23, 4}: \car \times_{\cb} \bigl( (G_1 \times^{B_{J_1}} L_{J_1}) \times G_2 \bigr) \rightarrow G_1 \times G_2.$$

We may define in the same way the variety $\car \times_{\cb} \bigl( G_1 \times (G_2 \times^{B_{J_2}} L_{J_2}) \bigr)$ and the proper morphism $f_{1, 2, 34}: \car \times_{\cb} \bigl( G_1 \times (G_2 \times^{B_{J_2}} L_{J_2}) \bigr) \rightarrow G_1 \times G_2$.

Now define an isomorphism $\iota: \car \times_{\cb} \bigl( (G_1 \times^{B_{J_1}} L_{J_1}) \times G_2 \bigr) \rightarrow \car \times_{\cb} \bigl( G_1 \times (G_2 \times^{B_{J_2}} L_{J_2}) \bigr)$ by $$\bigl((r_1, r_2), (g, g'), g_2 \bigr) \mapsto \Bigl( \bigl(r_1, r_2 ((g') \i, \th_{\d}(g') \i) \bigr), g, (g_2, \th_{\d}(g') \i) \Bigr)$$ for $r_1 \in \car_{\cc'}$ and $r_2 \in \car_{\cc}$. It is easy to see that $f_{1, 23, 4}=f_{1, 2, 34} \circ \iota$.

We have that \begin{align*} \pi_!(\bar{\mathbb Q}_l \odot \cl_{\mathbf w}) &=(f_{1, 23, 4})_! \bigl(\bar{\mathbb Q}_l \odot ((\cl_{w_1 s_i} \odot \cl'_{s_i}) \boxtimes \cl_{w_2}) \bigr) \\ &=(f_{1, 2, 34})_! \iota_! \bigl(\bar{\mathbb Q}_l \odot ((\cl_{w_1 s_i} \odot \cl'_{s_i}) \boxtimes \cl_{w_2}) \bigr) \\ &=(f_{1, 2, 34})_! \bigl(\bar{\mathbb Q}_l \odot (\cl_{w_1 s_i} \boxtimes (\cl_{w_2} \odot \cl'_{s_{\d(i)}})) \bigr).
\end{align*}

Now the lemma follows from \cite[4.2.1]{MS}.

\

Notice that for $\mathbf w \in W_1 \times W_2$, $$<\pi_! (\bar{\mathbb Q}_l \odot \cl_{\mathbf x})>_{\mathbf x \le \mathbf w}=<\pi_! (\bar{\mathbb Q}_l \odot \ca_{\mathbf x})>_{\mathbf x \le \mathbf w}.$$

Then we have the following consequence.

\begin{th28} Let $\mathbf w, \mathbf w' \in W_1 \times W_2$. Then 

(1) If $\mathbf w \rightarrow_{c, c'} \mathbf w'$ and $l(\mathbf w)>l(\mathbf w')$, then $$\pi_! (\bar{\mathbb Q}_l \odot \ca_{\mathbf w}) \in <\pi_! (\bar{\mathbb Q}_l \odot \ca_{\mathbf x})>_{l(\mathbf x)<l(\mathbf w)}.$$

(2) If $\mathbf w \approx_{c, c'} \mathbf w'$, then $$\pi_! (\bar{\mathbb Q}_l \odot \ca_{\mathbf w})=\pi_! (\bar{\mathbb Q}_l \odot \ca_{\mathbf w'}).$$
\end{th28}

\

Now combining the above results with Corollary 3.5, we have the following result which is a generalization of the key lemma in \cite[section 3]{H2}.

\begin{th29} Let $C$ be a unipotent character sheaf with respect to $\cc$ and $\cc'$. Then 

(1) $C \dashv \pi_! (\bar{\mathbb Q}_l \odot \cl_{\mathbf w})$ for some $\mathbf w$ that is of minimal length in the coset $W_{c'} \mathbf w W_c$.

(2) $C$ is a constitute of $\pi_! (\bar{\mathbb Q}_l \odot \ca_{\mathbf w})$ for some $\mathbf w$ that is of minimal length in the coset $W_{c'} \mathbf w W_c$.
\end{th29}

\

In the rest of this section, we consider the Hecke algebras.

\subsection*{6.6} Let $(W, I)$ be a Coxeter group. Given a map $L: I \rightarrow \mathbb Z$ with $L(i)=L(j)$ for all $i \neq j$ such that $m_{i j}$ is finite and odd. Let $\mathcal A=\mathbb Z[v, v \i]$, where $v$ is an indeterminate. Set $v_i=v^{L(i)} \in \mathcal A$.

Let $\mathcal H$ be the $\mathcal A$-algebra defined by the generators $T_{s_i}$ ($i \in I$) and the relations \begin{gather*}
\tag{a} (T_{s_i}-v_i)(T_{s_i}+v_i \i)=0 \text{ for } i \in I \\
\tag{b} T_{s_i} T_{s_j} T_{s_i} \cdots=T_{s_j} T_{s_i} T_{s_j} \cdots                                      \end{gather*} (both products have $m_{i j}$ factors) for any $i \neq j$ in $I$ such that $m_{i j}<\infty$. $\mathcal H$ is called the {\it Iwahori-Hecke algebra}.

For $w \in W$, we define $T_w=T_{s_{i_1}} T_{s_{i_2}} \cdots T_{s_{i_n}}$, where $w=s_{i_1} s_{i_2} \cdots s_{i_n}$ is a reduced expression.

For subset $J$ of $I$, we denote by $\ch_{J}$ the subalgebra of $\ch$ generated by $T_{s_j}$ ($j \in J$).

\subsection*{6,7} Now let $J, J' \subset I$ and $\d: W_J \rightarrow W_{J'}$ be an automorphism with $\d(J)=J'$. We assume furthermore that $L(j)=L(\d(j))$ for $i \in J$. Then there is a unique algebra isomorphism $D: \ch_{J_1} \rightarrow \ch_{J'_1}$ such that $D(T_{s_j})=T_{s_{\d(j)}}$ for $j \in J$.

\

Now we have the following result which is a generalization of some results in \cite{GP1} and \cite{GKP}.

\begin{th30} We keep the notation of the previous section. Let $\zeta: \ch \rightarrow \ca$ be a $\ca$-linear map such that $\zeta(h' h)=\zeta(h D(h'))$ for $h \in \ch$ and $h' \in \ch_J$. Let $\co$ be a $W_J$-orbit in $W$, where the $W_J$-action on $W$ is defined in Corollary 3.8. Let $w, w' \in \co_{\min}$. If moreover $W_J$ is a finite Coxeter group or $\co \cap W^{J'} \neq \varnothing$, then $\zeta(T_w)=\zeta(T_{w'})$.
\end{th30}

\begin{rmk} Some functions satisfying the condition in the Proposition arises in the study of parabolic character sheaves. See \cite[section 31]{L2}.
\end{rmk}

By Corollary 3.8, $w \sim_{\d} w'$. Now it suffices to prove the statement for $w'=x w \d(x) \i$ where $x \in W_J$ and either

(a) $l(x w)=l(x)+l(w)$;

or (b) $l(w \d(x) \i)=l(x)+l(w)$.

We only prove the case (a). Case (b) can be showed in the same way.

It is then easy to see that $T_x T_w=T_{x w}=T_{w' \d(x)}=T_{w'} T_{\d(x')}=T_{w'} D(T_x)$. Hence $\zeta(T_w)=\zeta(T_x \i (T_x T_w))=\zeta(T_x T_w D(T_x) \i)=\zeta(T_{w'})$. The proposition is proved.

\section{cuspidal $\s$-conjugacy classes}

In this section, we study the $\s$-conjugacy classes of finite Weyl group of type ABD. We will combine the approach in \cite[section 3]{GP2} and Corollary 3.8 to obtain a new way to understand the $\s$-conjugacy classes.

\subsection{} Let $\s: W \rightarrow W$ be an automorphism with $\s(I)=I$. For $w \in W$, set $\supp_{\s}(w)=\cup_{n \ge 0} \s^n \supp(w)$. Then $\supp_{\s}(w)$ is a $\s$-stable subset of $I$.

A $\s$-conjugacy class $\co$ of $W$ is called {\it cuspidal} if $\co \cap W_J=\varnothing$ for all proper $\s$-stable subset $J$ of $I$.

Let $V$ be the vector space spanned by $\a_i$ (for $i \in I$). We regard $W$ as a subgroup of $GL(V)$ and $\s$ as an element in $GL(V)$ in the natural way. For $w \in W$, set $$p_{w, \s}(q)=\det(q \cdot id_V-w \s).$$

Then it is easy to see that $p_{w, \s}(q)=p_{w', \s}(q)$ if $w$ is $\s$-conjugate to $w'$.

As in \cite[Exercise 1.15]{GP2}, for $w \in W$ and $i \in I$, define the length function $l_i(w)$ as the number of generators in $I$ conjugate to $s_i$ occuring in a reduced expression of $W$. Let $d$ be the minimal positive integer such that $\s^d(i)=i$. Set $$l_{i, \s}(w)=\sum_{k=0}^{d-1} l_{\s^k(i)} (w).$$

Then it is easy to see that if $w \approx_{\s} w'$, then $l_{i, \s}(w)=l_{i, \s}(w')$ for all $i \in I$.

\begin{th31}  If $W_J$ is finite for any proper $\s$-stable subset $J$ and $p_{w, \s}(1) \neq 0$, then the $\s$-conjugacy class of $w$ is cuspidal.
\end{th31}

\begin{rmk} This is a generalization of \cite[Lemma 3.1.10]{GP2}.
\end{rmk}

If $w \in W_J$ for some proper $\s$-stable subset $J$ of $I$, then $\supp_{\s}(w) \neq I$. Set $v=\sum_{i \notin \supp_{\s}(w)} \a_i$. Then $w \s(v)=w v=v+\a$ for some $\a \in \sum_{i \in \supp_{\s}(w)} \mathbb R \a_i$. Since $w \s$ is of finite order, we may assume that $(w \s)^n=id_V$. Thus $\sum_{1 \le i \le n} (w \s)^n v=n v+\b$ for some $\b \in \sum_{i \in \supp_{\s}(w)} \mathbb R \a_i$ and $\sum_{1 \le i \le n} (w \s)^n v$ is an eigenvector of $w \s$ with eigenvalue 1. Hence $p_{w, \s}(1)=0$.

\

The following Lemmas are obvious and we omit the proofs.

\begin{th32} Let $w \in W_J$ and $x \in W^J$, then $l(x w \s(x) \i) \ge l(v)$. In particular, if $J=\s(J)$ and $w$ is of minimal length in its $\s \mid_J$-conjugacy class in $W_J$, then it is of minimal length in its $\s$-conjugacy class in $W$.
\end{th32}

\begin{th33} Let $w, w' \in W$ with $w \rightarrow_{\s} w'$. Then $\supp_{\s}(w') \subset \supp_{\s}(w)$. In moreover, $w \approx_{\s} w'$, then $\supp_{\s}(w')=\supp_{\s}(w)$.
\end{th33}

\

Now we state the main theorem in this section which is a generalization of \cite[Theorem 3.2.7]{GP2}.

\begin{th34} For a finite Coxeter group $(W, I)$ and an automorphism $\s: W \rightarrow W$ with $\s(I)=I$, the following holds:

(P1) Let $w \in W$ be such that $\supp_{\s}(w)=I$ and that $\Cyc_{\s}(w)$ is terminal. Then the $\s$-conjugacy class of $w$ in $W$ is cuspidal and $w \in \co_{\min}$.

(P2) Let $\co$ be a cuspidal $\s$-conjugacy class of $W$. Then $\co_{\min}=\Cyc_{\s}(w)$ for any $w \in \co_{\min}$.

(P3) Let $\co, \co'$ be cuspidal $\s$-conjugacy classes of $W$ and $w \in \co_{\min}$, $w' \in \co'_{\min}$. Then $\co=\co'$ if and only if $p_{w, \s}(q)=p_{w', \s}(q)$ and $l_{i, \s}(w)=l_{i, \s}(w')$ for all $i \in I$.
\end{th34}

\begin{rmk} There exists cuspidal conjugacy classes $\co \neq \co'$ in finite Coxeter group of type $F_4$ such that $p_{w, \s}(q)=p_{w', \s}(q)$ for $w \in \co$ and $w' \in \co'$. See \cite[Appendix B]{GP2}.
\end{rmk}

\

As a consequence, it implies the Geck-Kim-Pfeiffer theorem \cite[2.6]{GKP}.

\begin{th35} Let $(W, I)$ be a finite Coxeter group and $\s: W \rightarrow W$ be an automorphism with $\s(I)=I$. Let $\co$ be a $\s$-conjugacy class of $W$. Then

(a) For each $w \in \co$, there exists an element $w' \in \co_{\min}$ such that $w \rightarrow_{\s} w'$.

(b) Let $w, v \in \co_{\min}$. Then there exists an element $w' \in \Cyc_{\s}(w)$ and an element $x \in W$ such that $w'$ is elementarily strongly $\s$-conjugate to $v$ via $x$. In particular, any two elements in $\co_{\min}$ are strongly $\s$-conjugate.
\end{th35}

\begin{rmk} The proof is similar to \cite[3.2.9]{GP2}.
\end{rmk}

By \cite[2.8]{GKP}, it suffices to prove the theorem for irreducible groups.

(a) Let $w \in \co$. If there exists $w' \in W$ such that $\supp_{\s}(w')$ is a proper subset of $I$ and $w \rightarrow_{\s} w'$. Then by induction on $\sharp I$, we have $w' \rightarrow_{\s} w''$ for some $w'' \in W_{\supp_{\s}(w')}$ which is of minimal length in its $\s$-conjugacy class in $W_{\supp_{\s}(w')}$. By Lemma 7.3, $w''$ also has minimal length in its $\s$-conjugacy class in $W$ and therefore $w \rightarrow_{\s} w'' \in \co_{\min}$.

Otherwise, $\supp_{\s}(w')=I$ for all $w' \in W$ with $w \rightarrow_{\s} w'$. Now let $w' \in W$ be such that $\Cyc_{\s}(w')$ is terminal and $w \rightarrow_{\s} w'$. Then by (P1) of the previous theorem, $\co$ is cuspidal and $w' \in \co_{\min}$.

Part (a) is proved.

(b) Since $w, v \in \co_{\min}$, we have that $l(w)=l(v)$ and $a w \s(a) \i=v$ for some $a \in W$. Write $a$ as $a=x b$ for $x \in W^{\supp_{\s}(w)}$ and $b \in W_{\supp_{\s}(w)}$. Set $w'=b w \s(b) \i$. Then $w' \in W_{\supp_{\s}(w)}$ and $v=x w' \s(x) \i$. By Lemma 7.3, $l(w') \le l(v)$. However, since $v \in \co_{\min}$, we have that $l(w')=l(v)$. Moreover, since $x \in W^{\supp_{\s}(w)}$, we have that $l(x w')=l(x)+l(w')$. Hence $w'$ is elementarily strongly $\s$-conjugate to $v$.

Since $w$ has minimal length in its $\s$-conjugacy class in $W_{\supp_{\s}(w)}$, its cyclic shift class $\Cyc_{\s}(w)$ is terminal. Hence by (P1) of the previous theorem, the $\s$-conjugacy class of $w$ in $W_{\supp_{\s}(w)}$ is cuspidal. Since $l(w')=l(w)$, by (P2) of the previous theorem, $w' \in \Cyc_{\s}(w)$.

Part (b) is proved.

\

Below is another generalization of the main theorem.

\begin{th36} Let $W$ be a finite Coxeter group and $\s$ be an automorphism of $W$ with $\s(I)=I$ and $\s^2=id$. Then for $w \in W$, $w$ and $\s(w) \i$ are in the same $\s$-conjugacy class.
\end{th36}

\begin{rmk} This is a generalization of \cite[Corollary 3.2.14]{GP2}. The proof is similar to {\it loc. cit.} and is omitted here.
\end{rmk}

\subsection*{7.8} We will prove the main theorem for Coxeter groups of classical type. The exceptional groups with $\s=id$ have been settled in \cite[Appendix B]{GP2} by direct computation. (P1) and (P2) of the main theorem have been settled for $^3 D_4$, $^2 F_4$ and $^2 E_6$ by direct computation in \cite[section 6]{GKP}. As to (P3), we can see from Table I-III in \cite[section 6]{GKP} that except for two classes in $^2 E_6$, minimal length elements in different cuspidal $\s$-conjugacy classes have different length. The only exception is the $\s$-conjugacy class of $w_1=s_1 s_3 s_1 s_2 s_4 s_3 s_1 s_5 s_4 s_3 s_1 s_6 s_5 s_4 s_3 s_1$ and the $\s$-conjugacy class of $w_2=s_2 s_4 s_5 s_4 s_2 s_3 s_4 s_5 s_6 s_5 s_4 s_2 s_3 s_4 s_5 s_6$. We have that $p_{w_1, \s}(q)=(q+1)^4 (q^2+q+1)$ and $p_{w_2, \s}(q)=(q^2+q+1)^2$. Thus (P3) also holds for these cases.

The Coxeter groups of classical type with $\s=id$ were first proved in \cite{GP1} and then in \cite{GP2} using cuspidal classes. We will give a new proof for these cases. We will also prove the main theorem for classical type with $\s \neq id$. 

The most difficult part of our proof is to find representatives of $$Cusp_{\s}(W)=\{w \in W; \supp_{\s}(w)=I, \Cyc_{\s}(w) \text{ is terminal }\}/\approx_{\s}.$$
We will find the representatives case by case. The general strategy is as follows.

Let $w \in W$ with $\supp_{\s}(w)=I$ and $\Cyc_{\s}(w)$ terminal. We
choose a maximal proper subset $J$ of $I$. Then $w \approx_{\s} w_1
v$ for some $w_1 \in W^{\s(J)}$ and $v \in W_{I(w_1, \s \mid_J)}$.
By Lemma 7.9 and 7.10 below, $\supp_{\s}(w_1)=I$ and the $\s
\Ad(w_1)$-conjugacy class of $v$ in $W_{I(w_1, \s \mid_J)}$ is
cuspidal. By induction on $I$, we may assume that $v$ is a
representative in $Cusp_{\s \Ad(w_1)} (W_{I(w_1, \s \mid_J)})$ that
we have found. In particular, $w \approx_{\s} w_1 v_1 v_2$ for $v_1
\in W^{\s w(K)}$ with $\supp_{\s \Ad(w_1)}(v_1)=I(w_1, \s\mid_J)$
and $v_2 \in W_{I(w_1 v_1, \s \mid K)}$. By Lemma 7.11 below, since
$\Cyc_{\s}(w)$ is terminal, $w_1$ and $v_1$ must satisfy some
condition.

In this way, we find some elements $x_k$ in $W$ such that for $w \in W$ with $\supp_{\s}(w)=I$ and $\Cyc_{\s}(w)$ terminal, we have $w \approx_{\s} x_k$ for some $x_k$. Now we calculate $p_{x_k, \s}(q)$ and check that

(1) $p_{x_k, \s}(1) \neq 0$;

(2) $p_{x_k, \s}(q) \neq p_{x_{k'}, \s}(q)$ for $x_k \neq x_{k'}$.

By 7.1 and Lemma 7.2, the $\s$-conjugacy class of $x_k$ is cuspidal and different $x_k$ belongs to different cuspidal class. Since each $\s$-conjugacy class contains at least one terminal cyclic shift class, $\Cyc_{\s}(x_k)$ is terminal for all $x_k$. Thus these $x_k$ are representatives of $Cusp_{\s}(W)$ and also representatives of cuspidal $\s$-conjugacy classes. (P1)-(P3) of the main theorem also hold in this case.

\begin{th37} Let $W$ be an irreducible Coxeter group.
Let $J \subset I$, $w \in W^{\s(J)}$ and $v \in W_{I(w, \s
\mid_J)}$. Then $\supp_{\s}(w v)=I$ if and only if
$\supp_{\s}(w)=I$.
\end{th37}

It is easy to see that if $\supp(w) \subset \supp(w v)$. Thus
$\supp_{\s}(w)=I$ implies that $\supp_{\s}(w v)=I$. On the other
hand, if $\supp_{\s}(w) \neq I$ and $\supp_{\s}(w v)=I$, then
$\cup_{n \ge 0} \s^n I(w, \s \mid_J) \supset I-\supp_{\s}(w)$. It is
easy to see that for $i \notin \supp_{\s}(w)$, $w \a_i$ is of the
form $\a_i+\sum_{j \in \supp(w)} a_j \a_j$ for some $a_j \in \mathbb
N \cup \{0\}$. By the definition of $I(w, \s \mid_J)$, we have that
$w \a_i=\a_i$ and $\s(i)=i$ for all $i \in I(w, \s
\mid_J)-\supp_{\s}(w)$. Therefore $I(w, \s \mid_J)-\supp_{\s}(w)$ is
$\s$-stable and $I(w, \s \mid_J) \supset I-\supp_{\s}(w)$.

Now since $W$ is irreducible, there exists $i \in I(w, \s \mid_J)-\supp_{\s}(w)$ and $j \in \supp_{\s}(w)$ such that $m_{i j} \neq 2$. Since $\s$ is an automorphism of $W$ and $\s(i)=i$, we have that $m_{i j}=m_{i, \s(j)}$. Therefore there exists $j \in \supp(w)$ such that $m_{i j} \neq 2$. Now let $w=s_{i_1} s_{i_2} \cdots s_{i_n}$ be a reduced expression and $m=\max\{k; m_{i, i_k} \neq 2\}$. Then $w
\a_i=s_{i_1} \cdots s_{i_m} \a_i=s_{i_1} \cdots s_{i_{m-1}} (\a_i+a \a_{i_m})$ for some $a \in \mathbb N$. Therefore $w \a_i=(s_{i_1} \cdots s_{i_{m-1}}) \a_i+a (s_{i_1} \cdots s_{i_{m-1}}) \a_{i_m}=\a_i+\sum_{j \in \supp(w)} a_j \a_j+a \a$, where $a_j \in \mathbb N \cup \{0\}$ and $\a=(s_{i_1} \cdots s_{i_{m-1}}) \a_{i_m}$ is a positive root. In particular, $w \a_i \neq \a_i$. That is a contradiction.

The lemma is proved.

\

Unless otherwise stated, we assume that $I=\{1, 2, \cdots, n\}$. Set $$s_{[a, b]}=\begin{cases} s_a s_{a-1} \cdots s_b, & \text{ if } a \ge b; \\ 0, & \text{ otherwise }. \end{cases}$$

\begin{th38} Let $(W, I)$ be an irreducible Coxeter group and $\s: W \rightarrow W$ be an automorphism with $\s(I)=I$. Let $d<n$ with

(a) if $\s^n(i) \le b-1$ for some $i \le b-1$, then $\s^n(i)=i$ for all $i \le b-1$;

(b) $m_{i j}=\d_{1, \mid i-j \mid}$ for $1 \le i, j \le b$;

(c) $m_{i, i'}=0$ for $i \le b-1$ and $i' \ge b+1$;

(d) $m_{b, \s^n(i)}=0$ for $i \le b-1$ and $n \in \mathbb Z$ with $\s^n(i) \neq i$.

Let $a \le b-1$. Let $w=\s \i(s_{[b, a]}) \i w_1 s_{[b, 1]} v_1 v_2$ with $w_1, v_1, v_2 \in W$, $\supp(w_1), \supp(v_1) \subset I-\cup_{n \in \mathbb N} \s^n \{1, 2, \cdots, b\}$, $\supp(v_2) \subset \{a+1, a+2, \cdots, b-1\}$ and $l(w)=2 b-a+1+l(w_1)+l(v_1)+l(v_2)$. Then $\Cyc_{\s}(w)$ is not terminal.
\end{th38}

\begin{rmk} Let $J=\{1, 2, \cdots, b-1\}$. The idea of the proof is to use the procedure in section 2 to obtain an element of the form $v_1 w_1$ where $w_1 \in {}^J W$ and $v_1 \in W_{I(w_1 \i, \s \mid_J)}$ such that $w \rightarrow_{\s \mid_J} v_1 w_1$ and $l(v_1 w_1)<l(w)$.
\end{rmk}

Set $x=s_{\s \i(b)} w_1 s_{[b, 1]} v_1$ and $y=v_2 s_{[b, a]} \i$. Then $x s_{\s^n(i)}=s_{\s^n(i)} x$ for $i \in \{1, 2, \cdots, b-1\}$ and $n \in \mathbb Z$ with $\s^n(i) \neq i$.

Notice that $x y x \i=s_{[b, 1]} y s_{[b, 1]} \i$. Then $l(x y x \i)=l(x)$. Now let $n_0$ be the minimal positive integer such that $\s^{n_0}(i)=i$ for $i \le b-1$. Then \begin{align*} w \rightarrow_{\s} & x y=(x y x \i) x \rightarrow_{\s} x \s(x y x \i)=\s(x y x \i) x \rightarrow_{\s} x \s^2(x y x \i) \\ \rightarrow_{\s} & \cdots \rightarrow_{\s} x \s^{n_0} (x y x \i)=x (x y x \i). \end{align*}

We can show in this way that $w \rightarrow_{\s} x (x^{a-1} y x^{-(a-1)})$. Notice that \begin{align*} x (x^{a-1} y x^{-(a-1)}) &=(x^a v_2 x^{-a}) x \bigl(x^{a-1} s_{[b-1, a]} \i x^{-(a-1)}) \bigr) \\ &=(x^a v_2 x^{-a}) x s_{[b-a, 1]} \i.
\end{align*}

Since $l((x^a v_2 x^{-a}) x s_{[b-a, 1]} \i) \le l(v_2)+l(x)-(b-a)<l(w)$, $\Cyc_{\s}(w)$ is not terminal. The lemma is proved.

\begin{th39} We keep the assumption in the previous lemma. Let $$w=\s \i(s_{[b, a]}) \i w_1 s_{[b, 1]} s_b v_1 s_{[b, a+1]} v_2$$ with $\supp(w_1), \supp(v_1), \supp(v_2) \in I-\cup_{n \in \mathbb N} \s^n \{1, 2, \cdots, b\}$ and $l(w)=l(w_1)+l(v_1)+l(v_2)+3 b-2 a+2$. If $2 a<b$, then $\Cyc_{\s}(w)$ is not terminal.
\end{th39}

\begin{rmk} This is a generalization of the ``Block exchange'' lemma in \cite[Lemma 3.4.5]{GP2}. The proof here is similar to the previous lemma.
\end{rmk}

Set $x=s_{\s \i(b)} w_1 s_{[b, 1]} s_b v_1$. Then \begin{align*} w & \rightarrow_{\s} x s_{[b, a+1]} v_2 s_{[b-1, a]} \i=x s_{[b, a+1]} s_{[b-1, a]} \i v_2=x s_{[b-2, a]} \i s_{[b, a]} v_2 \\ &=s_{[b-3, a-1]} \i x s_{[b, a]} v_2. \end{align*}

As in the proof of the previous lemma, we can show that $$s_{[b-3, a-1]} \i x s_{[b, a]} v_2 \rightarrow_{\s} x s_{[b, a]} v_2 s_{[b-3, a-1]} \i.$$ If $2 a<b$, then we can show in the same way that \begin{align*} w & \rightarrow_{\s} x s_{[b, 2]} v_2 s_{[b-2 a+1, 1]} \i=x s_{[b, 2]} s_{[b-2 a+1, 1]} \i v_2=x  s_{[b-2 a, 1]}  \i s_{[b, 1]} v_2.\end{align*}

Since $l(x s_{[b-2 a, 1]} \i s_{[b, 1]} v_2) \le l(v_2)+l(x)+b-(b-2 a)<l(w)$, $\Cyc_{\s}(w)$ is not terminal. The lemma is proved.

\

Now we will prove the main theorem for each type. We will use the
same labelling of Dynkin diagram as in \cite{Bo}. Set $J=I-\{\s
\i(1)\}$. Then $\s(J)=I-\{1\}$.

\

{\bf Type $A_n$}

\subsection*{7.12} Let $w \in W$ with $\supp_{id}(w)=I$ and $\Cyc_{id}(w)$ is terminal.
By Corollary 3.8, $w \approx_{id \mid_J} w_1 v$ for some $w_1 \in
W^J$ and $v \in W_{I(x, id \mid_J)}$. By Lemma 7.9,
$\supp_{id}(w_1)=I$. Thus $w_1=s_{[n, 1]}$. Then $I(x, id
\mid_J)=\varnothing$ and $w \approx_{id \mid_J} w_1$. It is easy to
see that $p_{w_1, id}(q)=\sum_{1 \le i \le n} q^i$. So there exists
a unique cuspidal conjugacy class, which is just the conjugacy class
that contains $s_{[n, 1]}$.

\

{\bf Type $^2 A_n$}

\begin{th40} Let $W$ be a Weyl group of type $A_n$ and $\s$ be an automorphism of order 2 on $W$ with $\s(I)=I$. For any sequence $\a=(\a_1, \a_2, \cdots, \a_l)$ with $\a_1 \ge \a_2 \ge \cdots \ge \a_l \ge 1$ and $\sum_{1 \le i \le l} (2 \a_i-1)=n+1$, we set $$w_{\a}=s_{[n+1-\a_1, 1]} s_{[n+2-\a_1-\a_2, \a_1+1]} \cdots s_{[n+l-\sum_{1 \le i \le l} \a_i, \sum_{1 \le i \le l-1} (\a_i)+1]}.$$ Let $w \in W$ with $\supp_{\s}(w)=I$ and $\Cyc_{\s}(w)$ is terminal. Then $w \approx_{\s} w_{\a}$ for some $\a$.
\end{th40}

We argue by induction on $n$. By Corollary 3.8, $w \approx_{\s \mid_J} w_1 v$
for some $w_1 \in W^{\s(J)}$ and $v \in W_{I(x, \s \mid_J)}$. By
Lemma 7.9, $\supp_{\s}(w_1)=I$. Thus $w_1=s_{[n+1-\a_1, 1]}$ for some
$\a_1 \ge 1$. Then $I(w_1, \s \mid_J)=\{\a_1+1, \a_1+2, \cdots,
n+1-\a_1\}$ and $\s \Ad(w_1)$ is an order-2 bijection on $I(w_1, \s
\mid_J)$. By Lemma 7.10, $v$ is contained in a cuspidal $\s
\Ad(w_1)$-conjugacy of $W_{I(w_1, \s \mid_J)}$. By induction
hypothesis, $v \approx_{\s \Ad(w_1)} s_{[n+2-\a_1-\a_2, \a_1+1]}
\cdots s_{[n+l-\sum_{1 \le i \le l} \a_i, \sum_{1 \le i \le l-1}
(\a_i)+1]}$ for some sequence $\a'=(\a_2, \a_3, \cdots, \a_l)$ with
$\a_2 \ge \a_3 \ge \cdots \ge \a_l \ge 1$ and $\sum_{2 \le i \le l}
(2 \a_i-1)=n+2-2 \a_1$. By Corollary 3.8,
$$w \approx_{\s \mid_J} s_{[n+1-\a_1, 1]} s_{[n+2-\a_1-\a_2,
\a_1+1]} \cdots s_{[n+l-\sum_{1 \le i \le l} \a_i, \sum_{1 \le i \le
l-1} (\a_i)+1]}.$$

Notice that $\Cyc_{\s}(w)$ is terminal. By Lemma 7.11, we have that $\a_1 \ge \a_2$. The lemma is proved.

\subsection*{7.14} We have that \begin{align*} p_{w_{\a}, \s}(q) &=\det(q \cdot id_V-w_{\a} \s)=\det(q \cdot id_V+w w_0) \\ &=(-1)^n \det(-q \cdot id_V-w w_0)=\frac{(-1)^n}{(-q-1)} \prod_{1 \le i \le l} ((-q)^{2 \a_i-1}-1) \\ &=\frac{1}{(q+1)} \prod_{1 \le i \le l} (q^{2 \a_1}+1).\end{align*} Thus $w_{\a}$ is contained in a cuspidal $\s$-conjugacy class. It is also easy to see that $p_{w_{\a}, \s}(q) \neq p_{w_{\a'}, \s}(q)$ for $\a \neq \a'$. By the argument in 7.8, the main theorem holds in this case.

We also showed that the cuspidal $\s$-conjugacy classes of $A_n$ are parametrized by the sequence $\a=(\a_1, \a_2, \cdots, \a_l)$ with $\a_1 \ge \a_2 \ge \cdots \ge \a_l \ge 1$ and $\sum_{1 \le i \le l} (2 \a_i-1)=n+1$. In other words, the cuspidal $\s$-conjugacy classes of $A_n$ are parametrized by the partitions of $n+1$ with only odd parts.

\

{\bf Type $B_n$}

\begin{th41} Let $W$ be a Weyl group of type $B_n$. For any sequence $\a=(\a_1, \a_2, \cdots, \a_l)$ with $\a_1 \ge \a_2 \ge \cdots \ge \a_l \ge 1$ and $\sum_{1 \le i \le l} \a_i=n$, we set $$w_{\a}=(s_{[n-1, \a_1]} \i s_{[n, 1]}) (s_{[n-1, \a_1+\a_2]} \i s_{[n, \a_1+1]}) \cdots (s_{[n, \sum_{1 \le i \le l-1} (\a_i)+1]}).$$ Let $w \in W$ with $\supp_{id}(w)=I$ and $\Cyc_{id}(w)$ is terminal. Then $w \approx_{id} w_{\a}$ for some $\a$.
\end{th41}

We argue by induction on $n$. By Corollary 3.8, $w \approx_{id \mid_J} w_1 v$ for some $w_1 \in W^J$ and $v \in W_{I(x, id \mid_J)}$. By Lemma 7.9, $\supp_{id}(w_1)=I$. Thus $w_1=s_{[n-1, \a_1]} \i s_{[n, 1]}$ for some $\a_1 \ge 1$. Then $I(w_1, id \mid_J)=\{\a_1+1, \a_1+2, \cdots, n\}$ and $\Ad(w_1)$ is the identity map on $I(w_1, id \mid_J)$. By Lemma 7.10, $v$ is contained in a cuspidal conjugacy of $W_{I(w_1, id \mid_J)}$. By induction hypothesis, $v \approx_{id} (s_{[n-1, \a_1+\a_2]} \i s_{[n, \a_1+1]}) \cdots (s_{[n, \sum_{1 \le i \le
l-1} (\a_i)]})$ for some sequence $\a'=(\a_2, \a_3, \cdots, \a_l)$ with $\a_2 \ge \a_3 \ge \cdots \ge \a_l \ge 1$ and $\sum_{2 \le i \le l} \a_i=n-\a_1$. By Corollary 3.8, $$w \approx_{J, id} (s_{[n-1, \a_1]} \i s_{[n, 1]}) (s_{[n-1, \a_1+\a_2]} \i s_{[n, \a_1+1]}) \cdots (s_{[n, \sum_{1 \le i \le l-1} (\a_i)+1]}).$$

Notice that $\Cyc_{\s}(w)$ is terminal. By Lemma 7.11, we have that $\a_1 \ge \a_2$. The lemma is proved.

\subsection*{7.16} By \cite[3.4.3]{GP2}, $p_{w_{\a}, id}(q)=\prod_{1 \le i \le l} (q^{\a_i}+1)$. Thus $w_{\a}$ is contained in a cuspidal conjugacy class. Moreover, $p_{w_{\a}, id}(q) \neq p_{w_{\a'}, id}(q)$ for $\a \neq \a'$. By the argument in 7.8, the main theorem holds in this case. We also showed that the cuspidal conjugacy classes of $B_n$ are parametrized by the partitions of $n$.

\

{\bf Type $^2 B_2$}

\subsection*{7.17} There is one cuspidal $\s$-conjugacy class, which is the class that contains $s_1 s_2 s_1$. The other minimal length element in the class is $s_2 s_1 s_2 \approx_{\s} s_1 s_2 s_1$. The main theorem holds in this case. We have that $p_{s_1 s_2 s_1, \s}(q)=(q+1)^2$.

\

{\bf Type $D_n$ and $^2 D_n$}

\subsection*{7,18} Let $0 \le a<b \le n$. Define $$w_{a, b}=\begin{cases} s_{[n-2, b]} \i s_{[n, a+1]}, & \text{ if } b \le n-1; \\ s_{[n-1, a+1]}, & \text{ if } b=n. \end{cases}$$

For any sequence $\a=(\a_1, \a_2, \cdots, \a_l)$ with $\a_1 \ge \a_2 \ge \cdots \ge \a_l \ge 1$ and $\sum_{1 \le i \le l} \a_i=n$, we set $$w'_{\a}=w_{0, \a_1} w_{\a_1, \a_1+\a_2} \cdots w_{\sum_{1 \le i \le l-1} \a_i, \sum_{1 \le i \le l} \a_i}.$$

\begin{th42} Let $W$ be a Weyl group of type $D_n$. Let $\s_0=id$ and $\s_1$ be the automorphism of order 2 on $W$ with $\s_1(I)=I$. Let $w \in W$ with $\supp_{\s_i}(w)=I$ and $\Cyc_{\s_i}(w)$ is terminal. Then $w \approx_{\s_i} w'_{\a}$ for some $\a$ with $2 \mid l-i$.
\end{th42}

We argue by induction on $n$. For $n=3$, it is easy to check that
the statement holds. Now assume that $n \ge 4$. By Corollary 3.8, $w
\approx_{\s_i \mid_J} w_1 v$ for some $w_1 \in W^{\s_i(J)}$ and $v
\in W_{I(x, \s_i \mid_J)}$. By Lemma 7.9, $\supp_{\s_i}(w_1)=I$. Thus
$w_1=w_{[0, \a_1]}$ for some $\a_1$ with $1 \le \a_1 \le n+i-1$.
Then $$I(w_1, \s_i \mid_J)=\begin{cases} \{\a_1+1, \a_1+2, \cdots,
n\}, & \text{ if } \a_1 \le n-2; \\ \varnothing, & \text{ if }
\a_1>n-2. \end{cases}$$ and $\s_i \Ad(w_1)$ is the bijection of
order $2-i$ on $I(w_1, \s_i \mid_J)$. By Lemma 7.10, $v$ is contained
in a $\s_{1-i}$-cuspidal conjugacy of $W_{I(w_1, \s_i \mid_J)}$. By
induction hypothesis, $v \approx_{\s_{1-i}} w_{\a_1, \a_1+\a_2}
\cdots w_{\sum_{1 \le i \le l-1} \a_i, \sum_{1 \le i \le l} \a_i}$
for some sequence $\a'=(\a_2, \a_3, \cdots, \a_l)$ with $\a_2 \ge
\cdots \ge \a_l \ge 1$, $\sum_{2 \le i \le l} \a_i=n-\a_1$
and $2 \mid (l-1)-(1-i)=l+i$. By Corollary 3.8, $$w \approx_{\s_i \mid_J} w_{0,
\a_1} w_{\a_1, \a_1+\a_2} \cdots w_{\sum_{1 \le i \le l-1} \a_i,
\sum_{1 \le i \le l} \a_i}.$$

Notice that $\Cyc_{\s_i}(w)$ is terminal. By Lemma 7.11, we have that $\a_1 \ge \a_2$. The lemma is proved.

\subsection*{7.20} We use $\tilde s_1, \tilde s_2, \cdots, \tilde s_n$ for the standard generators of the Weyl group of type $B_n$. By \cite[1.4.8]{GP2}, we may regard $W$ as a subgroup of a Weyl group of type $B_n$ via $s_i \rightarrow \tilde s_i$ for $i \le n-1$ and $s_n \rightarrow \tilde s_n \tilde s_{n-1} \tilde s_n$. For any partition $\a$ of $n$ with even numbers of parts, the element $w'_{\a}$ in $W$ is just the element $w_{\a}$ of the Weyl group of type $B_n$. Thus $p_{w'_{\a}, id}(q)=p_{w_{\a},
id}=\prod_{1 \le i \le l} (q^{\a_i}+1)$. Therefore $w'_{\a}$ is contained in a cuspidal conjugacy class of $W$. Moreover, $p_{w'_{\a}, id}(q) \neq p_{w'_{\a'}, id}(q)$ for $\a \neq \a'$. By the argument in 7.8, the main theorem holds in this case. We also showed that the cuspidal conjugacy classes of $D_n$ are parametrized by the partitions of $n$ with even numbers of parts.

It is easy to see that $\s_1=\tilde s_n$ as an element in $GL(V)$. Moreover, for any partition $\a$ of $n$ with odd numbers of parts, $\tilde s_n w'_{\a}=w_{\a}$ in the Weyl group of type $B_n$. Thus \begin{align*} p_{w'_{\a}, \s_1}(q) &=\det(q \cdot id_V-w'_{\a} \tilde s_n)=\det(q \cdot id_V-\tilde s_n w'_{\a}) \\ &=\det(q \cdot id_V-w_{\a})=\prod_{1 \le i \le l} (q^{\a_i}+1).\end{align*} Therefore $w'_{\a}$ is contained in a cuspidal $\s_1$-conjugacy class of $W$. Moreover, $p_{w'_{\a}, \s_1}(q) \neq p_{w'_{\a'},
\s_1}(q)$ for $\a \neq \a'$. By the argument in 7.8, the main theorem holds in this case. We also showed that the cuspidal $\s_1$-conjugacy classes of $D_n$ are parametrized by the partitions of $n$ with odd numbers of parts.

\

{\bf Type $^3 D_4$}

\begin{th43} Let $W$ be the Weyl group of type $D_4$ and $\s$ be an automorphism of $W$ with $\s(s_1)=s_3$, $\s(s_3)=s_4$ and $\s(s_4)=s_1$. Let $w \in W$ with $\supp_{\s}(w)=I$ and $\Cyc_{\s}(w)$ is terminal. Then $w \approx_{\s} w'$ for some $w' \in \{s_2 s_1, s_3 s_2 s_1 s_3, s_3 s_2 s_1 s_2 s_3 s_2, s_1 s_2 s_4 s_3 s_2 s_1 s_2 s_4\}$.
\end{th43}

By Corollary 3.8, $w \approx_{\s \mid_J} w_1 v$ for some $w_1 \in W^{\s(J)}$ and
$v \in W_{I(w_1, \s \mid_J)}$. By Lemma 7.9, $\supp_{\s}(w_1)=I$. Thus
$$w_1 \in \{s_2 s_1, s_{[3, 1]}, s_4 s_2 s_1, s_{[4, 1]}, s_2 s_{[4,
1]}, s_1 s_2 s_{[4, 1]}\}.$$ Moreover $$I(J, w_1, \s)=\begin{cases}
\varnothing, & \text{ if } w_1 \in \{s_2 s_1, s_4 s_2 s_1, s_{[4,
1]}\}, \\ \{2, 3\}, & \text{ if } w_1=s_{[3, 1]}, \\ \{4\}, & \text{
if } w_1=s_2 s_{[4, 1]}, \\ \{2, 4\}, & \text{ if } w_1=s_1 s_2
s_{[4, 1]}. \end{cases}$$

If $w_1 \in \{s_2 s_1, s_4 s_2 s_1, s_{[4, 1]}\}$, then $v=1$ and $w' \approx_{\s} w_1$. Notice that $s_4 s_2 s_1 \xrightarrow{s_4}_{\s} s_2$ and $s_{[4, 1]} \xrightarrow{s_4}_{\s} s_3 s_2$. Thus $w \approx_{\s} s_2 s_1$.

If $w_1=s_{[3, 1]}$, then $v \approx_{\Ad(w_1) \s} v_1$, where $v_1 \in \{1, s_3, s_2 s_3 s_2\}$. Thus $w \approx_{\s} w_1 v_1$. Notice that $s_{[3, 1]} \xrightarrow{s_3}_{\s} s_2 s_1 s_4 \xrightarrow{s_2}_{\s} s_1 s_4 s_2 \xrightarrow{s_4}_{\s} s_1 s_2 s_1 \xrightarrow{s_2}_{\s} s_1$. So $w \approx_{\s} s_3 s_2 s_1 s_3$ or $w \approx_{\s} s_3 s_2 s_1 s_2 s_3 s_2$.

If $w_1=s_2 s_{[4, 1]}$, then $v=1$ or $v=s_4$. Notice that $s_2 s_{[4, 1]} \xrightarrow{s_2}_{\s} s_1 s_{[4, 1]} \xrightarrow{s_4}_{\s} s_1 s_3 s_2$. We also have that $s_2 s_{[4, 1]} s_4 \xrightarrow{s_4}_{\s} s_4 s_2 s_4 s_3 s_2 s_4$ and $\s^{-1} (s_4 s_2 s_4 s_3 s_2 s_4)=s_3 s_2 s_1 s_2 s_3 s_2$. Thus $w \approx_{\s} s_2 s_{[4, 1]} s_4 \approx_{\s} s_3 s_2 s_1 s_2 s_3 s_2$.

If $w_1=s_1 s_2 s_{[4, 1]}$, then $v \approx_{\Ad(w_1) \s} v_1$, where $v_1 \in \{1, s_2, s_2 s_4\}$. Thus $w \approx_{\s} w_1 v_1$. Notice that $s_1 s_2 s_{[4, 1]} \xrightarrow{s_1}_{\s} s_2 s_{[4, 1]} s_3 \xrightarrow{s_4}_{\s} s_2 s_4 s_3 s_2$ and $s_1 s_2 s_{[4, 1]} s_2 \xrightarrow{s_1}_{\s} s_2 s_{[4, 1]} s_2 s_4 \xrightarrow{s_2}_{\s} s_{[4, 1]} s_2 s_4 s_2 \xrightarrow{s_3}_{\s} s_4 s_2 s_1 s_4 s_2$. Thus $w \approx_{\s} s_1 s_2 s_4 s_3 s_2 s_1 s_2 s_4$.

The lemma is proved.

\subsection*{7.22} We have that \begin{gather*}
p_{s_2 s_1, \s}(q)=q^4-q^2+1, p_{s_3 s_2 s_1 s_3, \s}(q)=(q^2-q+1)^2, \\ p_{s_3 s_2 s_1 s_2 s_3 s_2, \s}(q)=(q+1)^2 (q^2-q+1), \\ p_{s_1 s_2 s_4 s_3 s_2 s_1 s_2 s_4, \s}(q)=(q^2+q+1)^2.
\end{gather*}

Set $\cw=\{s_2 s_1, s_3 s_2 s_1 s_3, s_3 s_2 s_1 s_2 s_3 s_2, s_1 s_2 s_4 s_3 s_2 s_1 s_2 s_4\}$. Then $w$ is contained in a cuspidal $\s$-conjugacy class of $W$ for $w \in \cw$. Moreover, $p_{w, \s}(q) \neq p_{w', \s}(q)$ for $w \neq w' \in \cw$. By the argument in 7.8, the main theorem holds in this case.

\

In the rest of this section, we study the ``good'' elements.

\subsection*{7.23} Let $w \in W$, we call $d$ the $\s$-order of $w$ if $d$ is the minimal positive integer such that $w \s(w) \cdots \s^{d-1}(w)=1$ and $\s^d=1$.

Let $B^+$ be the braid monoid associated with $(W, I)$. Then there is a canonical injection $f: W \rightarrow B^+$ that identify the generators of $W$ with the generators of $B^+$ and $f(w_1 w_2)=f(w_1) f(w_2)$ if $w_1, w_2 \in W$ and $l(w_1 w_2)=l(w_1)+l(w_2)$. We will simply write $\underline w$ for $f(w)$.

Now the automorphism $\s$ extends to an automorphism of $B^+$ (which we denote by the same symbol).

We call an element $w \in W$ of $\s$-order $d$ a {\it good element} if there exists a sequence $I_1 \supset I_2 \supset \cdots \supset I_l$ of $I$ such that $$\underline w \s(\underline w) \cdots \s^{d-1}(\underline w)=\underline w_{I_1}^2 \underline w_{I_2}^2 \cdots \underline w_{I_l}^2 \text{ in } B^+.$$

The ``good'' elements for $\s=id$ were introduced in \cite{GM}. The above generalization appeared in \cite{GKP}.

\begin{th44} (1) Let $W$ be the Weyl group of type $A_n$ and $\s$ be an automorphism of order 2 with $\s(I)=I$. Then for any $a \le n$, $$\underline s_{[n+1-a, 1]} \s(\underline s_{[n+1-a, 1]}) \cdots \s^{2 a-2}(\underline s_{[n+1-a, 1]}) \underline w_{\{a+1, a+2, \cdots, n+1-a\}}=\underline w_I.$$

(2) Let $W$ be the Weyl group of type $B_n$. Then for any $a \le n$, $$(\underline s_{[n-1, a]} \i \underline s_{[n, 1]})^a \underline w_{\{a+1, a+2, \cdots, n\}}=\underline w_I.$$

(3) Let $W$ be the Weyl group of type $D_n$. Then for $a \le n-2$, $$(\underline s_{[n-2, a]} \i \underline s_{[n, 1]})^a \underline w_{\{a+1, a+2, \cdots, n\}}=\underline w_I.$$

(4) Let $W$ be the Weyl group of type $D_n$ and $\s$ be the automorphism of order 2 with $\s(I)=I$. Then \begin{gather*} (\underline s_{[n, 1]})^{n-1}=\underline w_I; \\ \underline s_{[n-1, 1]} \s(\underline s_{[n-1, 1]}) \cdots \s^{n-1}(\underline s_{[n-1, 1]})=\underline w_I. \end{gather*}
\end{th44}

We will prove part (1). The rest of the lemma can be showed in the same way.

By direct calculation, $$s_{[n+1-a, 1]} \cdots \s^{2 a-1}( s_{[n+1-a, 1]}) w_{\{a+1, a+2, \cdots, n+1-a\}} (\a_i)=-\a_{n+1-i}$$ for each simple root $\a_i$. Thus $$s_{[n+1-a, 1]} \s(s_{[n+1-a, 1]}) \cdots \s^{2 a-1}( s_{[n+1-a, 1]}) w_{\{a+1, a+2, \cdots, n+1-a\}}=w_I.$$ Moreover, $(2 a-1) l(s_{[n+1-a, 1]})+l(w_{\{a+1, a+2, \cdots, n+1-a\}})=n (n+1)/2$. Now part (1) follows from the definition of $f: W \rightarrow B^+$.

Part (1) is proved.

\begin{th45} (1) Let $W$ be the Weyl group of type $A_n$ and $\s$ be an automorphism of order 2 with $\s(I)=I$. Then for any $\a=(\a_1, \a_2, \cdots, \a_l)$ with $\a_1 \ge \a_2 \ge \cdots \ge \a_l \ge 1$ and $\sum_{1 \le i \le l} (2 \a_i-1)=n+1$, $$\underline w_{\a} \s(\underline w_{\a}) \cdots \s^{2d-1}(\underline w_{\a})=\underline w_{I_1}^{e_1} \underline w_{I_2}^{e_2-e_1} \cdots \underline w_{I_l}^{e_l-e_{l-1}},$$ where $d$ is the least common multiple of $2 \a_i-1$, $e_i=2d/(2 \a_i-1)$ and
$$I_i=\{\sum_{1 \le k \le i-1} \a_k-i+2, \sum_{1 \le k \le i-1} \a_k-i+3, \cdots, n-\sum_{1 \le k \le i-1} \a_k\}.$$

(2) Let $W$ be the Weyl group of type $B_n$. Then for any partition $\a=(\a_1, \a_2, \cdots, \a_l)$ of $n$, $$(\underline w_{\a})^{d}=\underline w_{I_1}^{e_1} \underline w_{I_2}^{e_2-e_1} \cdots \underline w_{I_l}^{e_l-e_{l-1}},$$ where $d$ is the least common multiple of $\a_i$, $e_i=d/a_i$ and $$I_i=\{\sum_{1 \le k \le i-1} \a_k+1, \sum_{1 \le k \le i-1} \a_k+2, \cdots, n\}.$$

(3) Let $W$ be the Weyl group of type $D_n$. Let $\s_0=id$ and $\s_1$ be the automorphism of order 2 on $W$ with $\s_1(I)=I$. Then for any partition $\a=(\a_1, \a_2, \cdots, \a_l)$ of $n$ with $l-i$ even, $$\underline w'_{\a} \s(\underline w'_{\a}) \cdots \s^{2d-1}(\underline w'_{\a})=\underline w_{I_1}^{e_1} \underline w_{I_2}^{e_2-e_1} \cdots \underline w_{I_l}^{e_l-e_{l-1}},$$ where $d$ is the least common multiple of $\a_i$, $e_i=2 d/a_i$ and $$I_i=\begin{cases} \{\sum\limits_{1 \le k \le i-1} \a_k+1, \sum\limits_{1 \le k \le i-1} \a_k+2, \cdots, n\},
& \text{ if } \sum\limits_{1 \le k \le i-1} \a_k \le n-2, \\ \varnothing, & \text{ otherwise.} \end{cases}$$
\end{th45}

\begin{rmk} Part (2) and (3) were proved in \cite[Proposition 4.3.11]{GP2} and part (1) was conjectured in \cite[5.6]{GKP}.
\end{rmk}

We will prove part (1). The rest of the lemma can be showed in the same way.

We argue by induction on $l$. Let $J=I-\{1\}$. Then $w_{\a}$ is of the form $w_1 v$, where $w_1=s_{[n+1-\a_1, 1]} \in W^J$ and $v \in W_{I(w_1, \s \mid_J)}$. It is easy to see that $$\underline w_{\a} \s(\underline w_{\a}) \cdots \s^{2d-1}(\underline w_{\a})=\underline w_1 \s(\underline w_1) \cdots \s^{2d-1}(\underline w_1) \underline v_1 \s_1(\underline v_1) \cdots \s_1^{2d-1}(\underline v_1)$$ where $\s_1=\Ad(w_1) \s$ and \begin{align*} v_1 &=\Ad \bigl(\s(w_1) \s^2(w_1) \cdots \s^{2d-1}(w_1) \bigr) \i v \\ &=\Ad(s_{[n,
\a_1]}) \Ad(s_{[n+1-\a_1, 1]} \i s_{[n, \a_1]})^{d-1} v=\Ad(s_{[n, \a_1]}) v.
\end{align*}

Notice that $\s_1$ is an order-2 automorphism on $W_{[\a_1, \a_1+1, \cdots, n-\a_1]}$.
By induction hypothesis, $$\underline v_1 \s_1(\underline v_1) \cdots \s_1^{2d'-1}(\underline v_1)=(\underline w_{I_2}^{e'_2} \underline w_{I_3}^{e'_3-e'_2} \cdots \underline w_{I_l}^{e'_l-e'_{l-1}}),$$ where $d'$ is the least common multiple of $2 \a_i-1$ for $i \ge 2$ and $e'_i=2 d'/(2 \a_i-1)$. By \cite[Proposition 4.1.9]{GP2}, \begin{align*} v_1 \s_1(v_1) \cdots \s_1^{d-1}(v_1) &=(\underline w_{I_2}^{e'_2})^{\frac{d}{d'}} (\underline w_{I_3}^{e'_3-e'_2})^{\frac{d}{d'}} \cdots
(\underline w_{I_l}^{e'_l-e'_{l-1}})^{\frac{d}{d'}} \\ &=\underline w_{I_2}^{e_2} \underline w_{I_3}^{e_3-e_2} \cdots \underline w_{I_l}^{e_l-e_{l-1}}. \end{align*}

By the previous lemma, \begin{align*} & \underline w_1 \s(\underline w_1) \cdots \s^{2d-1}(\underline w_1)
\underline w_{I_2}^{e_1} \\ &=\underline w_1 \s(\underline w_1) \cdots \s^{2d-e_1-1}(\underline w_1) \s \bigl(\underline w_1 \s(\underline w_1) \cdots \s^{e_1-1} \underline(w_1) \underline w_{\s(I_2)} \bigr) \underline w_{I_2}^{e_1-1} \\ &=\underline w_1 \s(\underline w_1) \cdots \s^{2d-e_1-1}(\underline w_1) \underline w_I \underline w_{I_2}^{e_1-1} \\ &=\underline w_1 \s(\underline w_1) \cdots \s^{2d-e_1-1}(\underline w_1) \underline w_{\s(I_2)}^{e_1-1} \underline w_I \\ &=\cdots=\underline w_I^{e_1}.
\end{align*}

Part (1) is proved.

\begin{th46} Let $(W, I)$ be a finite Coxeter group and $\s$ be an automorphism on $W$ with $\s(I)=I$. Let $\co$ be a $\s$-conjugacy class. Then there exists a good element $w \in \co_{\min}$.
\end{th46}

\begin{rmk} The non-twisted cases were proved in \cite{GM}. The twisted cases except the type $^2 A_n$ were proved in \cite{GKP}. The type $^2 A_n$ follows from the the previous corollary.
\end{rmk}

\bibliographystyle{amsalpha}

\end{document}